\documentclass{pnastwostrip}
\usepackage[pdftex]{graphicx}
\usepackage[T1]{fontenc}
\usepackage[utf8]{inputenc}
\usepackage{pnastwof}
\usepackage{amssymb,amsfonts,amsmath}
\usepackage{titlesec}
\usepackage[outline]{contour}
\usepackage{xcolor}
\usepackage{xargs}

\newcommand*{\cF}{\mathcal{F}}

\newcommand*{\cO}{\mathcal{O}}

\newcommand*{\bbC}{\mathbb{C}}

\newcommand*{\bbN}{\mathbb{N}}

\newcommand*{\bbR}{\mathbb{R}}

\newcommand*{\rL}{\mathrm{L}}

\def\I{\mathrm{i}}
\def\E{\mathrm{e}}
\newcommand*{\rth}{\mathrm{th}}
\def\imath{\mathrm{i}}

\newcommand{\be}[1]{\begin{equation} #1 \end{equation}}
\newcommand{\bes}[1]{\begin{equation*} #1 \end{equation*}}

\newcommand{\eas}[1]{\begin{align*} #1 \end{align*}}

\newcommand{\defn}[1]{\begin{definition} #1 \end{definition}}
\newenvironment{mybullet}{\begin{list}{$\bullet$ \quad}{}
	\setlength{\topsep}{0pt}
	\setlength{\parskip}{0pt}
	\setlength{\partopsep}{0pt}
	\setlength{\parsep}{0pt}         
	\setlength{\itemsep}{0pt}}{\end{list}}

\newenvironment{myenumerate}{\begin{enumerate}
	\setlength{\topsep}{0pt}
	\setlength{\parskip}{0pt}
	\setlength{\partopsep}{0pt}
	\setlength{\parsep}{0pt}         
	\setlength{\itemsep}{0pt}}{\end{enumerate}}

\renewenvironment{theorem}[1][Theorem]{\addtocounter{theorem}{1}\begin{trivlist}
	\item[\hskip \labelsep {\it\bfseries #1 \thetheorem.}] \it}{\end{trivlist}}

\newenvironment{theorem*}[1][Theorem]{\begin{trivlist}
	\item[\hskip \labelsep {\it\bfseries #1.}] \it}{\end{trivlist}}

\newcommand*{\spn}{\mathrm{span}}

\newcommand*{\ord}[1]{\cO \left ( #1 \right ) }

\newcommand{\res}[1]{\ifmmode{#1\!\times\!#1}\else{#1$\!\!\;\times\!\!\;$#1}\fi}

\makeatletter
\newcommand{\eg}{e.g.\@ifnextchar\ {}{\ }}
\newcommand{\ie}{i.e.\@ifnextchar\ {}{\ }}
\newcommand{\etc}{etc.\@ifnextchar\ {}{\ }}
\makeatother

\renewcommand*{\eqref}[1]{[Eq.\:\ref{#1}]}

\def\eref#1{\eqref{#1}}
\def\R#1{\eqref{#1}}
\def\fref#1{Fig.\:\ref{#1}}
\def\freftwo#1#2{Figs.~\ref{#1} and~\ref{#2}}

\makeatletter
\def\@cite#1#2{{(#1\if@tempswa,#2\fi)}}
\makeatother

\newcommand{\ifnonempty}[2]{\ifx\empty#1\empty\else#2\fi}
\newcommand{\ifempty}[2]{\ifx\empty#1\empty#2\else\fi}

\newcommand{\fontsf}[2][7]{\font\@fontsf\frutiger at #1pt\@fontsf#2}
\newcommand{\fontsfb}[2][7]{\font\@fontsfb\frutigerbold at #1pt\@fontsfb#2}

\def\colpt#1{c@{\hspace{#1pt}}}
\def\colper#1{c@{\hspace{#1\linewidth}}}

\makeatletter

\newlength{\@colwidth}
\newlength{\@imgstroke}
\newlength{\@overfontsize}
\font\@overfont\frutigerbold at \@overfontsize
\newlength{\@capfontsize}
\newlength{\@caplinesize}
\font\@capfont\frutiger at \@capfontsize
\newlength{\@sidefontsize}
\newlength{\@sidelinesize}
\font\@sidefont\frutigerbold at \@sidefontsize
\def\@overcolor{white}
\def\@overstrokecolor{black}
\contourlength{1pt}

\def\colwidth#1{\setlength{\@colwidth}{#1}}

\def\imgstroke#1{\setlength{\@imgstroke}{#1}}
\def\overfontsize#1{\setlength{\@overfontsize}{#1}}
\def\capfontsize#1{\setlength{\@capfontsize}{#1}}
\def\caplinesize#1{\setlength{\@caplinesize}{#1}}
\def\sidefontsize#1{\setlength{\@sidefontsize}{#1}}
\def\sidelinesize#1{\setlength{\@sidefontsize}{#1}}
\def\capfont#1{\font\@capfont#1 at \@capfontsize}
\def\overfont#1{\font\@overfont#1 at \@overfontsize}
\def\sidefont#1{\font\@sidefont#1 at \@sidefontsize}
\def\overcolor#1{\def\@overcolor{#1}}
\def\overstrokecolor#1{\def\@overstrokecolor{#1}}

\newcommand{\captext}[1]{\fontsize{\@capfontsize}{\@caplinesize}\selectfont\@capfont#1}
\newcommand{\sidetext}[1]{\fontsize{\@sidefontsize}{\@sidelinesize}\selectfont\@sidefont#1}

\newcommand{\vertsidetext}[2]{\raisebox{#1}{\rotatebox[origin=c]{90}{\sidetext{#2}}}}

\newcommand{\img}[2]{\setlength{\fboxsep}{0pt}\setlength{\fboxrule}{\@imgstroke}%
	\fbox{\includegraphics[width=#1\linewidth]{#2}}}
\newcommand{\overimg}[3]{\setlength{\fboxsep}{0pt}\setlength{\fboxrule}{\@imgstroke}%
	\fbox{\begin{overpic}#1\put(#2){#3}\end{overpic}}}

\newcommand{\overtext}[3][r]{%
	\setlength{\unitlength}{\@colwidth}%
	\begin{picture}(1,1)%
	\put(0,0){\setlength{\fboxsep}{0pt}\setlength{\fboxrule}{\@imgstroke}%
		\fbox{\includegraphics[width=\@colwidth]{#3}}}%
	\put(0,0.03){%
		\makebox[\@colwidth][#1]{%
			\hspace{0.03\unitlength}%
			\color{\@overcolor}\contour{\@overstrokecolor}{\@overfont#2}%
			\hspace{0.03\unitlength}%
		}%
	}%
	\end{picture}%
}

\newcommandx{\overerr}[3][1={},2=r]{\overtext[#2]{\ifnonempty{#1}{Err 
#1\%}}{#3}}

\makeatother

\def\AAA#1{#1}
\def\BBB#1{#1}
\def\CCC#1{#1}

\titleformat{\parstripagraph}[runin]{\bfseries}{\theparagraph}{1em}{}
\titlespacing*{\paragraph}{0pt}{0.5ex}{1em}
\titleformat{\subsection}[runin]{\subsectionfont}{\thesubsection}{1em plus 
0.25em minus 0.25em}{}
\titlespacing*{\subsection}{0pt}{2ex plus 0.5ex minus 0.2ex}{1em}

\frutigerbold at 7.5pt

\makeatletter

\renewcommand*{\eqref}[1]{[Eq.\:\ref{#1}]}
\def\@cite#1#2{{(#1\if@tempswa,#2\fi)}}
\newcommand{\listpara}[1]{%
\paragraph{#1}%
}
\makeatother

\DeclareGraphicsExtensions{.png,.jpg,.pdf}
\graphicspath{
	{fig/res-dep/}
	{fig/tem/}
	{fig/fm/}
	{fig/fm-studer/}
	{fig/fixed-p-pom/}
	{fig/sk/}
	{fig/no-rip/}
	{fig/tv-shuffle-glpu/}
	{fig/tv-shuffle-pom/}
	{fig/incoh-rand/}
	{fig/Coherence/}
	{fig/fliptests/cells/}
	{fig/fliptests/geoeye/}
	{fig/fliptests/glpu/}
	{fig/fliptests/glpu_radial/}
	{fig/fliptests/radioi/}
	{fig/fliptests/watch/}
	{fig/fixed-n-glpu/}
	{fig/fixed-n-glpu-noise/}
	{fig/fixed-n-watch/}
	{fig/fixed-p-cells/}
	{fig/fixed-p-cells-bern/}
	{fig/fixed-p-glpu/}
	{fig/fixed-p-glpu-bern/}
	{fig/universality/}
	{fig/frames-watch/}
	{fig/modelcs/}
	{fig/structure/}
	{fig/}
	{fig/bogdan/}
}

\url{}
\copyrightyear{}
\issuedate{Issue Date}
\volume{Volume}
\issuenumber{Issue Number}

\begin{document}

\title{On asymptotic structure in compressed sensing}

\author{%
	Bogdan Roman\affil{1}{Department of Applied Mathematics and Theoretical 
	Physics, University of Cambridge, Cambridge CB3 0WA, UK},
	Ben Adcock\affil{2}{Department of Mathematics, Purdue University, West 
	Lafayette, IN 47907, USA}%
	\and%
	Anders Hansen\affil{1}{}%
}

\contributor{Submitted to Proceedings of the National Academy of Sciences
of the United States of America}

\maketitle

\begin{article}
\begin{abstract}
This paper demonstrates how new principles of compressed sensing, namely 
asymptotic incoherence, asymptotic sparsity and multilevel sampling, can be 
utilised to better understand underlying phenomena in practical compressed 
sensing and improve results in real-world applications. The contribution of 
the paper is 
fourfold: 
First, it explains how the sampling strategy depends not only on the signal 
sparsity but also on its structure, and shows how to design effective sampling 
strategies utilising this. 
Second, it demonstrates that the optimal sampling strategy and the efficiency 
of compressed sensing also depends on the resolution of the problem, and 
shows how this phenomenon markedly affects compressed sensing results and how 
to exploit it.
Third, as the new framework also fits analog (infinite dimensional) models 
that govern many inverse problems in practice, the paper describes how it 
can be used to yield substantial improvements.
Fourth, by using multilevel sampling, which exploits the structure of the 
signal, the paper explains how one can outperform random 
Gaussian/Bernoulli sampling even when the classical $l^1$ recovery algorithm 
is replaced by modified algorithms which aim to exploit structure such as 
model based or Bayesian compressed sensing or approximate message passaging. 
This final observation raises the question whether universality is desirable 
even when such matrices are applicable.
Examples of practical applications investigated in this paper include Magnetic 
Resonance Imaging (MRI), Electron Microscopy (EM), Compressive Imaging (CI) 
and Fluorescence Microscopy (FM). For the latter, a new compressed sensing 
approach is also presented.
\end{abstract}


\dropcap{C}ompressed sensing (CS), introduced by Cand\`{e}s, Romberg \& 
Tao~\cite{CandesRombergTao} and Donoho~\cite{donohoCS}, states that under 
appropriate conditions one can overcome the Nyquist sampling barrier and 
recover signals using far fewer samples than dictated by the classical Shannon 
theory. This has important implications in many practical 
applications which caused CS to be intensely researched in the past decade.

CS problems can be divided into two types. 
\textbf{Type~I} are problems where the physical device imposes the 
sampling operator, but allows some limited freedom to design the sampling 
strategy. This category is vast, with examples including Magnetic Resonance 
Imaging (MRI), Electron Microscopy (EM), Computerized Tomography, Seismic 
Tomography {and} Radio Interferometry.
\textbf{Type~II} are problems where the sensing mechanism offers freedom 
to design both the sampling operator and the strategy. Examples include 
Fluorescence Microscopy (FM) and Compressive Imaging (CI) (e.g.\ single pixel 
and lensless cameras). In these two examples, many practical setups 
still impose some restrictions regarding the sampling operator, e.g.\ 
measurements must typically be binary.

Traditional CS is based on three pillars: \textit{sparsity} (there are $s$ 
important coefficients in the vector to be recovered, however, the location is 
arbitrary), \textit{incoherence} (the values in the measurements matrix should 
be 
uniformly spread out) and \textit{uniform random subsampling}.

For Type~I problems the issue is that the above pillars are {often lacking}. 
As we will argue, {many} Type~I problems are coherent due to the 
physics or because they are infinite-dimensional. The traditional CS framework 
is simply not applicable. However, CS was used successfully in many Type~I
problems, though with very different sampling techniques than uniform random 
subsampling, which lack a mathematical justification.

For Type~II problems the traditional CS theory is applicable, e.g.\ in 
CI one can use random Bernoulli matrices. The issue is that the use of 
complete randomness does not allow one to exploit the structure of the signal 
to be recovered from a sampling point of view. As we argue, real world signals 
are not sparse, but asymptotically sparse in frames such as wavelets or their 
*-let cousins. In particular, the asymptotic sparsity is highly structured.

\paragraph{New CS principles.} To bridge the gap between theory and 
practice, 
the authors 
introduced a new CS theory~\cite{AHPRBreaking} that replaces the traditional 
CS pillars with three new CS principles: 
\textit{asymptotic incoherence}, \textit{asymptotic sparsity} and 
\textit{multilevel sampling}. The new theory and principles reveal that
the optimal sampling strategy and benefits of CS depend on two factors: the 
structure of the signal and the resolution. This suggests a new understanding 
of the underlying phenomena and of how to improve CS results in practical 
applications, which are main topics of this paper. At the same time, this 
paper demonstrates how the new CS principles go hand in hand even with 
applications where traditional CS is applicable, and that substantial 
improvements and flexibility can be obtained.

\section{Traditional Compressed Sensing}
\label{s:TradCS}

A traditional CS setup is as follows. The aim is to recover a signal $f$ from 
an incomplete (subsampled) set of measurements $y$. Here, $f$ is represented 
as a vector in $\bbC^N$ and is assumed to be $s$-sparse in some orthonormal 
basis $\Phi\in\bbC^{N\times N}$ (e.g.\ wavelets) called \textit{sparsity}
basis. This means that its vector of coefficients $x = \Phi f$ has at most 
$s$ nonzero entries. Let $\Psi\in\bbC^{N\times N}$ be an orthonormal basis, 
called \textit{sensing} or \textit{sampling} basis, and write 
$U=\Psi\Phi^*=(u_{ij})$, which is an isometry. The coherence of $U$ is
\be{\label{e:incoh}
\mu(U) = \max_{i,j} | u_{ij} |^2 \in[1/N,1].
}
and $U$ is said to be perfectly incoherent if $\mu(U) = 1/N$.

Let the \textit{subsampling pattern} be the set $\Omega \subseteq \{ 
1,\ldots,N \}$ of cardinality $m$ with its elements chosen uniformly at random.
Owing to a result by Cand\`{e}s \& Plan~\cite{Candes_Plan} and Adcock \& 
Hansen~\cite{BAACHGSCS}, if we have access to the subset of measurements $y = 
P_\Omega \Psi f$ then $f$ can be recovered from $y$ exactly with probability 
at least $1-\epsilon$ if
\be{
\label{e:cond_tradCS}
m \gtrsim \mu(U) \cdot N \cdot s \cdot \left(1+\log (1/\epsilon)\right) 
\cdot \log(N),
}
where $P_\Omega\in{\{0,1\}}^{N \times N}$ is the diagonal projection matrix 
with 
the $j^{\rth}$ entry $1$ if $j \in \Omega$ and 0 otherwise, and the notation 
$a\gtrsim b$ means that $a \geq C\, b$ where $C>0$ is some 
constant independent of $a$ and $b$. Then, $f$ is recovered by solving
\be{\label{e:l1}
\min_{z\in\bbC^N} \|z\|_1 \quad\text{subject to}\quad \|y-P_\Omega Uz\| \leq 
\eta.
}
where $\eta$ is chosen according to the noise level (0 if noiseless). The key 
estimate \R{e:cond_tradCS} shows that the number of measurements $m$ required 
is, up to a log factor, on 
the order of the sparsity $s$, provided the coherence $\mu(U) = 
\ord{1/N}$.  This is the case, for example, when $U$ is the DFT, which was 
studied in some of the first CS papers~\cite{CandesRombergTao}.


\begin{figure}
\centering
\colwidth{0.35\linewidth}
\begin{tabular}{@{}c@{~}c@{}}
\overtext[l]{(a)}{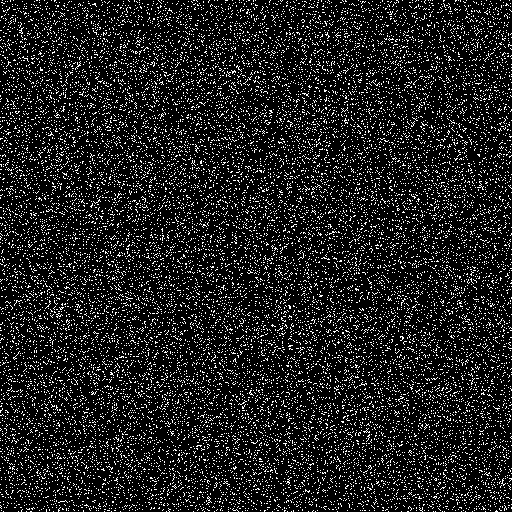}&
\overtext[l]{(b)}{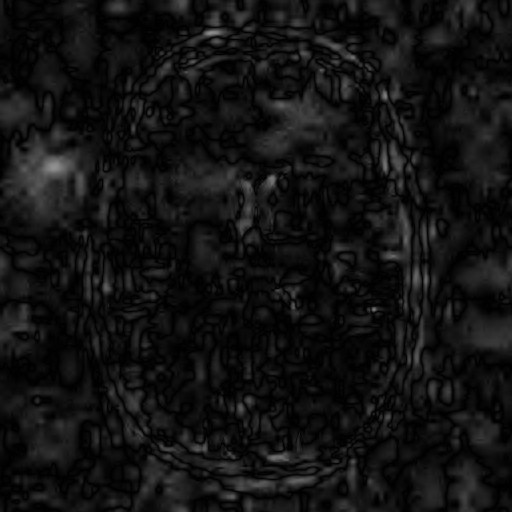}\\
\overtext[l]{(c)}{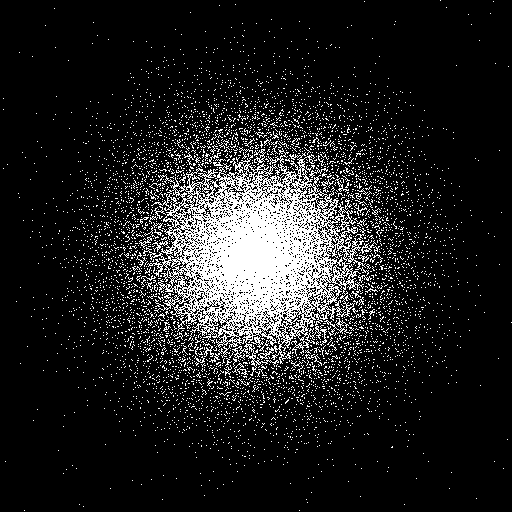}&
\overtext[l]{(d)}{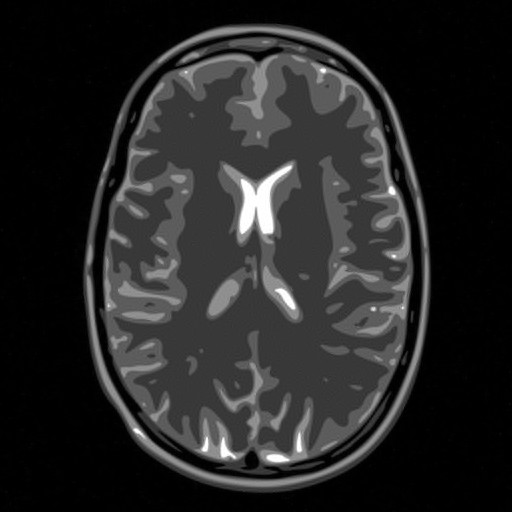}
\end{tabular}
\caption{(a) 12.5\% uniform random subsampling scheme, (b) CS 
reconstruction from uniform subsampling, (c) 12.5\% multilevel subsampling 
scheme, (d) CS reconstruction from multilevel subsampling.\label{f:CS_Evol}}
\end{figure}

\subsection{The real world is often coherent.}
\label{s:trad-incoh}

Consider the MRI CS setup, 
i.e.\ $U = \Psi_\mathrm{dft} \Phi^*_\mathrm{dwt} 
\in \mathbb{C}^{N \times N}$, where $\Psi_\mathrm{dft}$ and 
$\Phi_\mathrm{dwt}$ are the discrete Fourier and wavelet transforms 
(sampling and sparsity bases) respectively. The coherence here is
$$
\mu(U) = \ord{1}, \quad N \rightarrow \infty,
$$
for any wavelet basis, so this problem 
has the worst possible coherence. The traditional CS bound 
\eref{e:cond_tradCS} states that all samples are 
needed in this case (i.e.\ full sampling, $m=N$), even though the orignal 
signal is typically highly sparse in wavelets. This lack of incoherence means 
that uniform random subsampling leads to a very poor recovery.  This is 
known in MRI and is illustrated in \fref{f:CS_Evol}.

The root cause of this lack of incoherence is the discretization of what is 
intrinsically an infinite-dimensional problem into a finite-dimensional one. 
In short, $U$ {converges to an infinite matrix}~\cite{AHPRBreaking} and 
since the incoherence is the supremum of its entries, there exists some $N$ 
for which a \textit{coherence barrier} is hit, resulting in the worst case for 
a CS recovery. This is not restricted to MRI.
Any such discretization of an infinite-dimensional problem 
will suffer the same fate, including MRI, tomography, microscopy, seismology, 
radio interferometry etc. Changing $\Psi$ may provide marginal benefits, if 
any, since the coherence barrier always occurs at some $N$.

\begin{figure}[b]
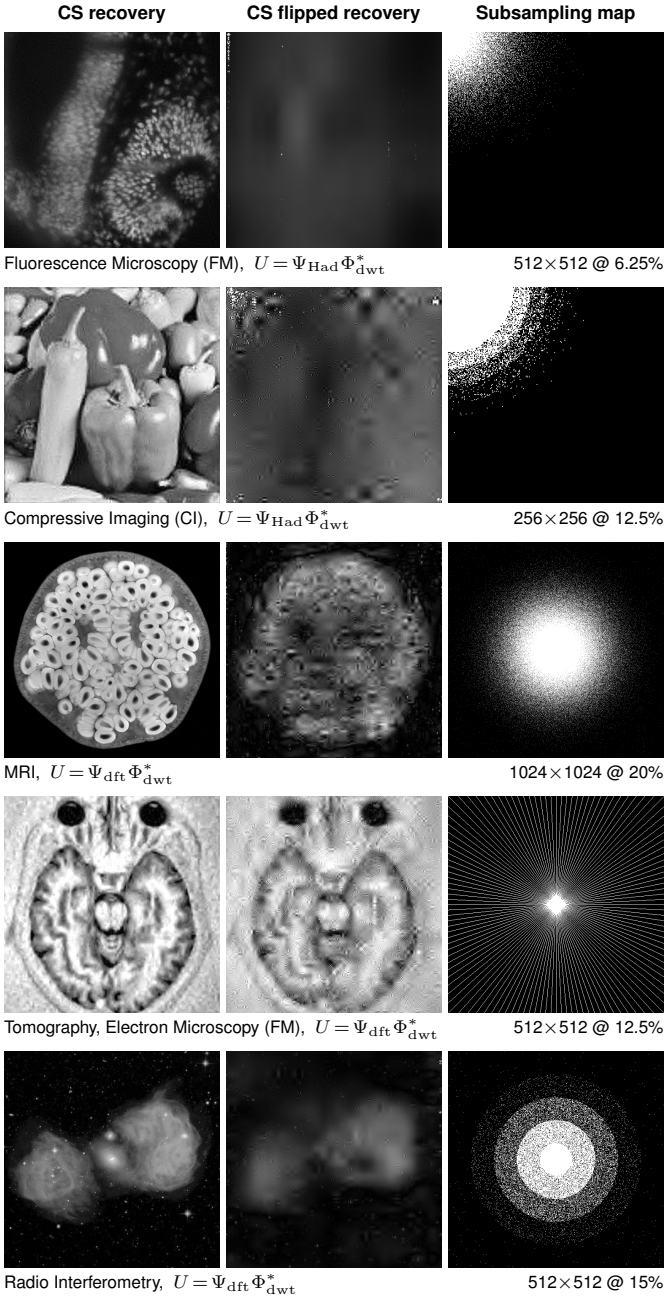

\centering
\frutiger at 6.5pt
\colwidth{0.33\linewidth}
\newcommand{\abrdesc}[5]{\\[-2pt]\multicolumn{3}{@{}l@{}}{%
	\fontsize{7pt}{7.2pt}\selectfont\abrdescfont%
	#5,~~$U\!=\!\Psi_\mathrm{#3}\Phi^*_\mathrm{#4}$%
	\hfill#1$\times$#1 @ #2\%}}%
\begin{tabular}{@{}\colper{0.01}\colper{0.01}c@{}}
\sidetext{CS recovery} & \sidetext{CS flipped recovery} & 
\sidetext{Subsampling map}\\[2pt]
  \overerr{0512_6p3_had_rwtdwtDB4_9_T8EP_rec}&
  \overerr{0512_6p3_had_rwtdwtDB4_9_flip1_WAD1_rec}&
  \overerr{0512_6p3_had_rwtdwtDB4_9_T8EP_map}
	\abrdesc{512}{6.25}{Had}{dwt}{Fluorescence Microscopy (FM)}\\[4pt]
  \overerr{0256_12p6_had_dwtdwtDB4_5_CEUF_rec}&
  \overerr{0256_12p6_had_wltdwtDB4_0_flip1_3XIM_rec}&
  \overerr{0256_12p6_had_dwtdwtDB4_5_CEUF_map}
	\abrdesc{256}{12.5}{Had}{dwt}{Compressive Imaging (CI)}\\[4pt]
  \overerr{1024_20p_dft_rwtdwtDB4_7_LAT2_rec}&
  \overerr{1024_20p_dft_rwtdwtDB4_7_flip1_NFY2_rec}&
  \overerr{1024_20p_dft_rwtdwtDB4_7_LAT2_map}
	\abrdesc{1024}{20}{dft}{dwt}{MRI}\\[4pt]
  \overerr{0512_12p_dft_rwtdwtDB4_6_L5GV_rec}&
  \overerr{0512_12p_dft_rwtdwtDB4_6_flip1_KTIG_rec}&
  \overerr{0512_12p_dft_rwtdwtDB4_6_L5GV_map}
	\abrdesc{512}{12.5}{dft}{dwt}{Tomography, Electron Microscopy (FM)}\\[4pt]
  \overerr{0512_15p_dft_db4_9_U9DD_rec.png}
  &\overerr{0512_15p_dft_db4_9_flip_B3NV_rec.png}
  &\overerr{0512_15p_dft_db4_9_U9DD_map.png}
	\abrdesc{512}{15}{dft}{dwt}{Radio Interferometry}
\end{tabular}
\caption{{\capfb Flip test.} Recovery from direct versus flipped wavelet coefficients 
showing that the RIP does not hold in these cases. The percentage shown is the 
subsampled fraction of Fourier/Hadamard coefficients.
\label{f:flip-tests}}
\end{figure}

\subsection{Sparsity, flip test and the absence of RIP.}
\label{s:trad-sparsity}

Traditional CS states that the sampling strategy is completely independent of 
the location of the nonzero coefficients of an $s$-sparse vector $x$, i.e.\ 
with the $s$ nonzero coefficients at arbitrary locations. 
The \textit{flip test} allows one to evaluate whether this holds in practice. 
Let 
$x\in\bbC^N$ be a {vector}, and $U\in\bbC^{N \times N}$ a measurement matrix.
We then sample according to some pattern $\Omega\subseteq\{1,\hdots,N\}$
with $|\Omega| = m$ and solve \R{e:l1} for $x$, \ie
$\min \|z\|_1$ s.t $P_{\Omega}Uz = P_{\Omega}Ux$ to obtain a 
reconstruction $z=\alpha$. Now we flip $x$ to obtain a vector $x'$ with 
reverse entries, $x'_i = x_{N-i}, i=1,\hdots,N$ and solve \R{e:l1} for $x'$ 
using the same $U$ and $\Omega$, i.e.\ $\min\|z\|_1$ s.t. 
$P_{\Omega}Uz = P_{\Omega}Ux'$. Assuming $z$ to be a solution, then by 
flipping $z$ we obtain a second reconstruction $\alpha'$ of the original 
vector 
$x$, where $\alpha'_i=z_{N-i}$.

Assume $\Omega$ is a sampling pattern for recovering $x$ using 
$\alpha$. If sparsity alone dictates the reconstruction quality, then 
$\alpha'$ must yield the same reconstruction quality (since $x'$ has the same 
sparsity as $x$, being merely a permutation of $x$). Is this 
true in practice? 

\fref{f:flip-tests} investigates this for several 
applications using $U = \Psi_\mathrm{dft}\Phi_\mathrm{dwt}^*$ 
or $U = \Psi_\mathrm{Had}\Phi_\mathrm{dwt}^*$, where $\Psi_\mathrm{dft}$, 
$\Psi_\mathrm{Had}$, $\Phi_\mathrm{dwt}$ are the discrete Hadamard, Fourier 
and wavelet transforms respectively. As is evident, the flipped recovery 
$\alpha'$ is substantially worse than its unflipped version $\alpha$. {This 
confirms that sparsity alone does not dictate the reconstruction quality.  
Furthermore, note that $P_{\Omega} U$ cannot satisfy an RIP for realistic 
values of $N$, $m$ and $s$.  Had this been the case, both vectors would have 
been recovered with the same error, and this is in direct contradiction with 
the results of the flip test.}

%
%
%

\begin{figure}
\colwidth{0.28\linewidth}
\imgstroke{0.25pt}
\centering
\begin{tabular}{@{}\colper{0.01}\colper{0.01}\colper{0}}
\overerr{0256_12p6_dft_tv_0p05_VLWV_map}&
\overerr{0256_12p6_dft_tv_0p05_VLWV_ref}&
\overerr[36.9]{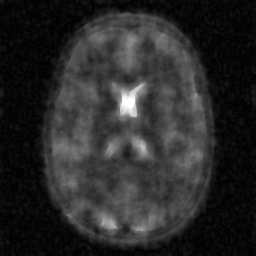}\\[-2pt]
\captext{Subsample map} & \captext{Original} & \captext{TV recovery}\\[4pt]
\overerr{0256_12p6_dft_tv_0p05_equiv_white_DKV7_map}&
\overerr{0256_12p6_dft_tv_0p05_equiv_white_DKV7_ref}&
\overerr[2.52]{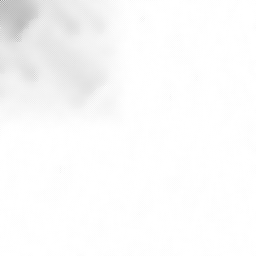}\\[-2pt]
\captext{Same subsample map} & \captext{Permuted gradients} & \captext{TV 
recovery}
\end{tabular}
\caption{{\capfb TV flip test.} TV recovery at \res{256} from 8192 DFT samples  
(12.5\% subsampling). The \textit{Permuted gradients} image 
was built from the gradient vectors of the \textit{Original} image,
having the same TV norm and gradient sparsity, differing only in the 
ordering and sign of the gradient vectors.
The large error difference confirms that sparsity structure matters for TV 
recovery as well.
\label{f:TV_shuffle}}
\end{figure}

It is worth noting that the same 
phenomenon exists for total variation (TV) minimization. 
Briefly, the CS TV problem is $\min_{z \in \bbC^n} \| z \|_{TV}$ s.t. $\| y - 
P_\Omega \Psi z \| \leq \eta$, where the TV norm $\|x\|_{TV}$ in case of 
images is the $\ell^1$ norm of the 
image gradient, $\|x\|_{TV} = \sum_{i,j} \|\nabla x(i,j)\|_2$ with $
\nabla x(i,j) = \{D_1x(i,j), D_2x(i,j)\}$, $D_1x(i,j) = x(i+1,j)-x(i,j)$, 
$D_2x(i,j) = x(i,j+1)-x(i,j)$. \fref{f:TV_shuffle} shows an experiment where we
we chose an image $x\in[0,1]^{N\times N}$ and then built an image 
$x'$ from the gradient of $x$ so that $\{\|\nabla x'(i,j)\|_2\}$ is a 
permutation of $\{\|\nabla x(i,j)\|_2\}$ for which $x'\in[0,1]^{N\times N}$.
Thus, the two images have the same ``TV sparsity'' and the same TV norm. It is 
evident how the reconstruction errors differ substantially for the 
two images when using the same sampling pattern, confirming that sparsity 
structure matters for TV recovery as well.

\section{New Compressed Sensing Principles}
\label{s:new}

The previous discussion on traditional CS calls for a more general
approach. We consider the generalization of the traditional principles of 
sparsity, incoherence, uniform random subsampling into \textit{asymptotic 
sparsity}, 
\textit{asymptotic incoherence} and \textit{multilevel 
subsampling}~\cite{AHPRBreaking}.

\subsection{Asymptotic sparsity.}
\label{s:new-sparsity}

We saw that signal structure is essential, but what structure describes 
such sparse signals? Let us consider a wavelet basis $\{ \varphi_n \}_{n 
\in \bbN}$.  Recall that there exists a decomposition of $\bbN$ into 
finite subsets according to different wavelet scales, 
i.e.\ $
\bbN = \bigcup_{k \in \bbN} \mathcal{M}_k,
$
where $\mathcal{M}_k=\{ M_{k-1}+1,\ldots,M_k \}$ is the 
set of indices corresponding to the $k^{\rth}$ scale, with $0 = M_0 < M_1 < 
M_2 < \ldots$.
Let $x \in l^2(\bbN)$ be the coefficients of a function $f$ in this basis.
Suppose that $\epsilon \in(0,1]$ is 
given, and define the \textit{local sparsities}
\begin{equation}\label{e:sk_def1}
s_k=s_k(\epsilon) = \min\Big\{L:\Big\|\!\sum_{i\in\mathcal{M}_{k,L}} \!\!
x_i\varphi_i\Big\| \geq 
\epsilon\, \Big\| \sum_{i\in \mathcal{M}_k} x_i\varphi_i  \Big\|\,\Big\},
\end{equation}
where $\mathcal{M}_{k,L}\subseteq \mathcal{M}_k$ is the set of indices of the 
largest $L$ coefficients at the $k^\rth$ scale, i.e.\ $|x_l|\geq|x_j|, \forall 
l\in \mathcal{M}_{k,L}, \forall j\in \mathcal{M}_k\setminus \mathcal{M}_{k,L}$.
In order words, $s_k$ is the effective sparsity of the wavelet 
coefficients of $f$ at the $k^{\rth}$ scale. 

\begin{figure}
\centering
\raisebox{20pt}{\includegraphics[width=0.3\linewidth]{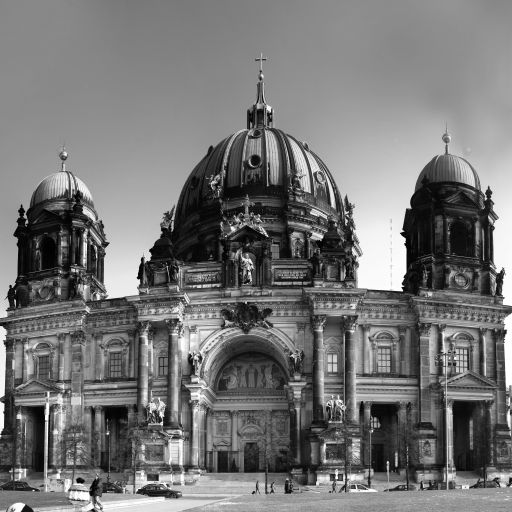}}%
\hspace{15pt}%
\includegraphics[width=0.6\linewidth]{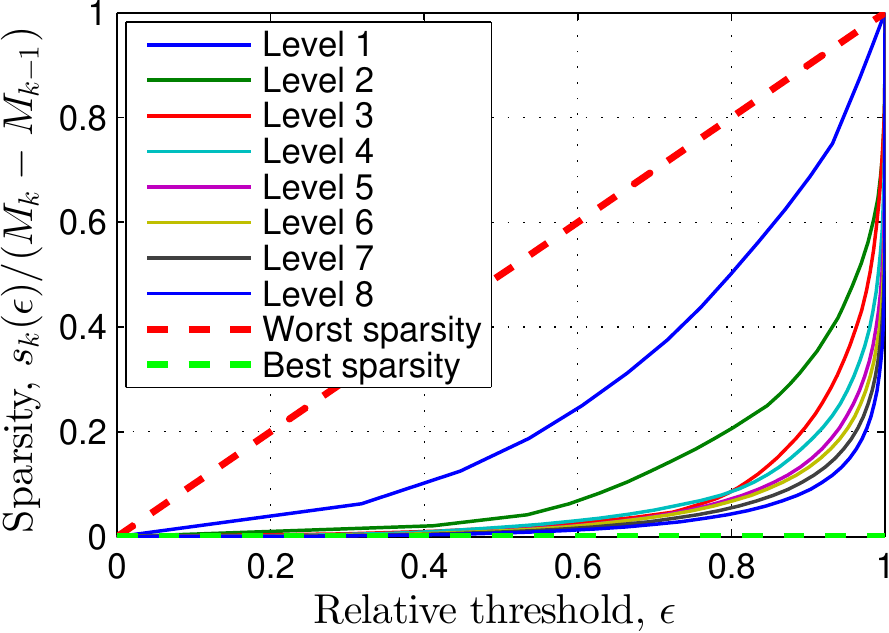}\\[4pt]
\caption{Sparsity of Daubechies-8 coefficients of an 
image. The levels correspond to wavelet scales and $s_k(\epsilon)$ is given by 
\R{e:sk_def1}. Each curve shows the relative sparsity at level $k$ as 
a function of $\epsilon$.  The decreasing nature of the curves for increasing 
$k$ confirms asymptotic sparsity \R{e:sk_decay}.\label{f:sk}}
\end{figure}

Sparsity of $x$ means that for a given large scale $r 
\in\bbN$, the ratio $s / M_r \ll 1$, where $M = M_r$ and 
$s = \sum_{k=1}^r s_k$ is the total sparsity of $x$. However, \fref{f:sk} 
shows that besides being sparse, practical signals have more structure,
namely \textit{asymptotic sparsity}, \ie
\be{\label{e:sk_decay}
s_k (\epsilon) / (M_k-M_{k-1}) \rightarrow 0,
}
rapidly as $k\!\rightarrow\!\infty$, $\forall\epsilon\!\in\!(0,1]$: they 
are far sparser at fine scales (large $k$) than at coarse scales (small 
$k$). This also holds for other function systems such as curvelets 
\cite{Cand}, contourlets~\cite{Vetterli} or shearlets 
\cite{Gitta}.

Given the structure of modern function systems such as wavelets and 
their generalizations, we propose the notion of sparsity in levels:
\defn{
\label{def:asy_sparsity}
Let $x\in\bbC$. For $r \in \bbN$ let $\mathbf{M} = (M_1,\ldots,M_r) \in 
\bbN^r$ and $\mathbf{s} = (s_1,\ldots,s_r) \in \bbN^r$, with $s_k \leq 
M_k - M_{k-1}$, 
$k=1,\ldots,r$, where $M_0 = 0$.  We say that $x$ is 
$(\mathbf{s},\mathbf{M})$-sparse if, for each $k=1,\ldots,r$, the sparsity band
\bes{
\Delta_k:=\mathrm{supp}(x) \cap \{ M_{k-1}+1,\ldots,M_{k} \},
}
satisfies $| \Delta_k | \leq s_k$.
We denote the set of $(\mathbf{s},\mathbf{M})$-sparse vectors by 
$\Sigma_{\mathbf{s},\mathbf{M}}$.
}

\subsection{Asymptotic incoherence.}
\label{s:new-incoh}

\begin{figure}
\fontsize{7pt}{8pt}\selectfont
\centering
\begin{tabular}{@{}\colper{0.025}\colper{0.01}\colper{0.01}\colper{0}}
\includegraphics[height=0.29\linewidth]{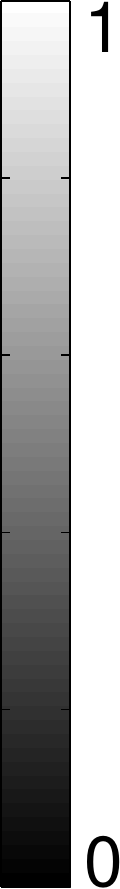}&
\includegraphics[width=0.29\linewidth]{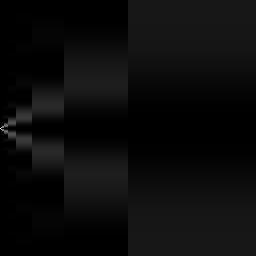}&
\includegraphics[width=0.29\linewidth]{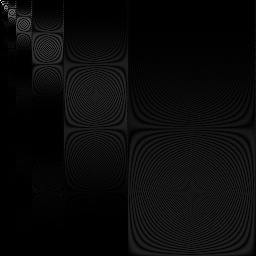}&
\includegraphics[width=0.29\linewidth]{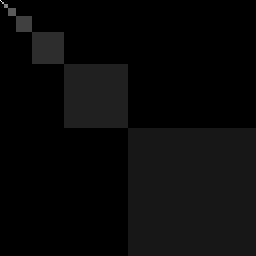}\\
&
{$U=\Psi_\mathrm{dft}\Phi^*_\mathrm{db3}$} &
{$U=\Psi_\mathrm{dct}\Phi^*_\mathrm{db3}$} &
{$U=\Psi_\mathrm{Had}\Phi^*_\mathrm{haar}$}
\end{tabular}
\caption{Visualizing incoherence. The absolute values of the matrix $U$. 
\label{f:incoh}}
\end{figure}

In contrast with random matrices (e.g.\ 
Gaussian or Bernoulli), many sampling and sparsifying operators typically 
found in practice yield fully coherent problems, such as the 
Fourier with wavelets case discussed earlier. \fref{f:incoh} shows the 
absolute values of the entries of 
the matrix $U$ for three examples. Although there are large values of $U$ in 
all three case (since $U$ is coherent as per \R{e:incoh}), these are 
isolated to a leading submatrix. Values get asymptotically smaller once we 
move away from this region. 

\defn{
Let $\{U_N\}$  be a sequence of isometries 
with $U_N \in \bbC^N$. $\{U_N\}$ is asymptotically incoherent if
$
\mu(P^{\perp}_K U_N),\ \mu(U_N P^{\perp}_K)\!\rightarrow\!0, 
$
when $K\!\rightarrow\!\infty$, with $N/K\!=\!c, \forall c\geq 1$.
Here $P_K$ is the projection onto $\spn \{e_j : j = 1,\ldots,K \}$, where $\{ 
e_j 
\}$ is the $\bbC^N$ canonical basis, and 
$P^{\perp}_K$ is its orthogonal complement.
}

In brief, $U$ is asymptotically incoherent if the coherences of the 
matrices formed by removing either the first $K$ rows or columns of $U$ are 
small. As \fref{f:incoh} shows, Fourier/wavelets, discrete cosine/wavelets and 
Hadamard/wavelets are examples of asymptotically incoherent problems.

\subsection{Multilevel sampling.}
\label{s:new-multilevel}

Asymptotic incoherence calls for a different strategy than uniform 
random sampling. High coherence in the first few rows of $U$ means that 
important information about the signal to be recovered 
is likely to be contained in the corresponding measurements, and thus we 
should fully sample these rows. Once outside this 
region, as coherence starts decreasing, we can subsample gradually.

%

\defn{
Let $r\!\in\!\bbN$, $\mathbf{N}\!=\!(N_1,\ldots,N_r)\!\in\!\bbN^r$ with 
$1\!\leq\!N_1\!<\!\ldots\!<\!N_r$, 
$\mathbf{m}\!=\!(m_1,\ldots,m_r)\!\in\!\bbN^r$, 
with $m_k\!\leq\!N_k-N_{k-1}$, $k=1,\ldots,r$, and suppose that
$\Omega_k\!\subseteq\!\{ N_{k-1}\!+\!1,\ldots,N_{k}\}, |\Omega_k|\!=\!m_k$,
are chosen uniformly at random, where $N_0\!=\!0$.  We refer to the set
$\Omega\!=\!\Omega_{\mathbf{N},\mathbf{m}}\!=\!\bigcup_{k=1}^r\Omega_k$
as an $(\mathbf{N},\mathbf{m})$-multilevel sampling scheme (using $r$ levels).
\label{def:multi_level}
}

Briefly, for a vector $x$, the sampling amount $m_k$ needed in 
each sampling band $\Omega_{k}$ is determined by the sparsity of $x$ in the 
corresponding sparsity band $\Delta_k$ and the asymptotic coherence 
$\mu(P^{\perp}_{N_k} U)$.

\section{Designing a multilevel sampling scheme}

Let $U$ be an isometry.
The $(k,l)^{\rth}$ local coherence of $U$ with respect to 
$\mathbf{N}$ and $\mathbf{M}$ is given by
\eas{
\mu_{\mathbf{N},\mathbf{M}}(k,l) &= 
\sqrt{\mu(P^{N_{k-1}}_{N_{k}}UP^{M_{l-1}}_{M_{l}}) \cdot  
\mu(P^{N_{k-1}}_{N_{k}}U)},
}
where $k,l=1,\ldots,r$ and $P^{a}_{b}$ is the projection matrix 
corresponding to indices $\{a+1,\hdots,b\}$.
The $k^{\rth}$ relative sparsity is
\bes{
S_k = S_k(\mathbf{N},\mathbf{M},\mathbf{s}) =  \max_{z \in 
\Theta}\|P_{N_k}^{N_{k-1}}U z\|^2,
}
where $\Theta = \{z : \|z\|_\infty \leq 1, 
|\mathrm{supp}(P_{M_l}^{M_{l-1}}z)| = s_l, \, l=1,\hdots, r\}$.

\paragraph{Sampling in levels.}
From~\cite{AHPRBreaking} we know that to recover an 
$(\mathbf{s},\!\mathbf{M})$-sparse vector $x\!\in\!\mathbb{C}^{N}$ from 
multilevel sampled measurements $P_{\Omega}Ux$,
{it suffices that the number of samples $m_k$ in each level} satisfies
\be{
\label{e:cond_levels}
1 \gtrsim \frac{N_k-N_{k-1}}{m_k} \cdot\!
\left(\sum_{l=1}^r \mu_{\mathbf{N},\mathbf{M}}(k,l) \cdot s_l\right) \!\cdot 
\log(N),
 }
 where
$
m_k \gtrsim \hat m_k \cdot \log(N),
$
and $\hat{m}_k$ is such that
\be{
\label{e:cond_hatm}
 1 \gtrsim \sum_{k=1}^r \left(\frac{N_k-N_{k-1}}{\hat m_k}- 1\right) \cdot 
 \mu_{\mathbf{N},\mathbf{M}}(k,l)\cdot \tilde s_k,
 }
{$\forall l=1,\ldots,r$ and $\forall\tilde{s}_1,\ldots,\tilde{s}_r\in 
(0,\infty)$ for which $\tilde s_{1}+ \hdots + \tilde s_{r}  \leq s_1+ \hdots + 
s_r$ and $\tilde s_k \leq S_k(\mathbf{N},\mathbf{M},\mathbf{s})$.}



The bounds \R{e:cond_levels} and 
\R{e:cond_hatm} are key. As opposed to the traditional CS bound 
\R{e:cond_tradCS}, 
which relates the total amount of subsampling $m$ to the global coherence and 
the global sparsity, these new bounds relate the local sampling amounts
$m_k$ to the local coherences $\mu_{\mathbf{N},\mathbf{M}}(k,l)$ and local and 
relative sparsities $s_k$ and $S_k$. A direct result is that the theorem 
agrees with the flip test shown earlier: the optimal sampling 
strategy must indeed depend on the signal structure. Another important note is 
that the bounds \R{e:cond_levels} and \R{e:cond_hatm} are sharp {in the sense 
that}
they reduce to the information-theoretic limits {in a number of important 
cases}.
Furthermore, in the case of Fourier sampling with wavelet sparsity, they 
provide near-optimal recovery guarantees using the infinite-dimensional 
generalization of the theorem. For further details and proofs, the 
reader is referred to~\cite{AHPRBreaking}.

We devised a flexible \textit{all-round} multilevel sampling scheme.
Assuming the coefficients $f\in\mathbb{C}^{N\times N}$ of a sampling 
orthobasis, e.g.\ DFT, our 
multilevel sampling scheme divides $f$ into $n$ regions delimited by 
$n-1$ equispaced concentric circles plus the full square. Normalizing the 
support of $f$ to $[-1,1]^2$, the circles have radius $r_k$ with 
$k=0,\ldots,n-1$, which are 
given by $r_0=m$ and $r_{k}=k\cdot\frac{1-m}{n-1}$ for $k>0$, where 
$m\in(0,1)$ is a parameter. In each of the $n$ regions, the fraction of 
coefficients sampled with uniform probability is
\begin{equation}
p_k=\exp\!\big(\!-\!(b\,k\!\!\;/\!\!\;n)^a\big),
\label{e:sublaw-gengauss}
\end{equation}
where $k=0,\ldots,n$ and $a>0$ and $b>0$ are parameters.
The total 
fraction of subsampled coefficients is $p=\sum_k p_k A_k$, where $A_k$ is the 
normalized area of the $k$th region. Since $p_0=1$ and $p_k>p_{k+1}$, the 
first region will sample all coefficients and the remaining regions will 
sample a decreasing fraction. An example is shown later in 
\ref{f:res-dep-submap}.

\section{Effects and benefits of the new principles}
\label{s:effects}

Having reviewed the theory, we now discuss the 
important effects and benefits of asymptotic incoherence and 
asymptotic sparsity, and of exploiting them via multilevel sampling. We show 
how they allow one to improve the CS recovery, and take practical examples 
from 
FM, MRI, EM and CI. We begin here with a short summary of 
the effects and benefits, which are detailed in subsequent sections.

  \listpara{The optimal sampling strategy is signal dependent.}
As the flip test shows, the optimal sampling strategy depends on the 
structure of the signal. Multilevel sampling 
takes this into account and allows one to further improve the CS recovery 
by tailoring the sampling according to e.g.\ the resolution and expected 
wavelet 
structure of the signal. This has additional advantages
as one can mitigate {application-specific} hurdles or target 
{application-specific} features (e.g.\ brain and spine imaging would use 
different subsampling schemes). This applies to both Type~I and Type~II 
problems. 
The \AAA{\textit{FM} and \textit{Resolution dependency}} sections below
provide examples from FM and MRI.

  \listpara{Resolution dependency.}
An important effect is that regardless of the sampling basis and 
subsampling scheme, the quality of the reconstruction increases as 
resolution increases. This is first revealed by fixing the subsampling 
strategy and fraction across resolutions, and a more {striking} result is 
obtained by fixing the number (instead of fraction) of samples, revealing 
hidden details, previously inaccessible. This is due to signals being 
typically 
increasingly (asymptotically) sparse at higher resolutions.
The \AAA{\textit{FM} and \textit{Resolution dependency}} sections
show this phenomena with examples from FM and MRI.

\listpara{Infinite dimensional CS.}
The new theory provides a good fit to {some} real-world problems 
that are fundamentally continuous, e.g.\ EM or MRI. 
The errors arising from recovering the continuous samples 
using discrete models are {sometimes} significant {\cite{GLPU}}. 
The section \AAA{\textit{Infinite dimensional problems}}
 discusses this aspect and 
shows an EM example where such noticeable 
errors can be overcome by using generalized sampling theory 
\cite{AHPRBreaking} and recovery into boundary wavelets.

\listpara{Structured sampling vs Structured Recovery.}
We exploit sparsity structure by using multilevel sampling of 
asymptotically incoherent matrices and standard $\ell^1$ minimization 
algorithms. Alternatively, sparsity structure can be exploited by using 
universal sampling matrices (e.g.\ random 
Gaussian/Bernoulli) and modified recovery algorithms. The section 
\AAA{\textit{Structured sampling vs Structured Recovery}}
discusses and compares the two approaches, highlighting the advantages of the 
former, {which, in contrast with the latter, allows to choose the 
sparsity frame, is applicable to both Type~I and Type~II problems, and 
yields overall superior results.}

\listpara{Structure vs Universality: Asymptotic vs Uniform 
incoherence.} 
The universality {property} of random sensing matrices (e.g.\ Gaussian, 
Bernoulli), explained later on, is a reason for their popularity in 
traditional CS. But is universality desirable when the signal {sparsity} 
is structured? {Should one use universal matrices when 
there is freedom to choose the sampling operator, i.e.\ in Type~II problems? 
Random matrices are largely inapplicable in Type~I problems, where the 
sampling operator is imposed. 
The \AAA{\textit{Structure vs Universality}} section argues
that 
universal matrices offer little room to exploit extra structure the signal may 
have, even in Type~II problems, and that non-universal matrices, such as 
Hadamard, coupled with multilevel sampling {provide} a better solution for 
both 
Type~I and Type~II problems as they can exploit the prevalent asymptotic 
sparsity of signals in practice.
}

\listpara{Storage and speed.}
Random matrices, popular in traditional CS, besides being 
inapplicable in many applications, are also slow and 
require (large) storage. This yields slow recovery and limits the maximum 
signal size{, which severely affects computations and, more importantly, 
sparsity structure.}
The \AAA{\textit{Storage/speed}} section
discusses this 
aspect and also shows that simply addressing the speed and storage 
problems via fast 
transforms and non-random matrices is not sufficient to achieve improved 
recovery compared to what multilevel sampling of non-universal matrices can 
offer.

\listpara{Frames and TV.}
Although investigating frames or TV as sparsity systems in CS is not the 
purpose of this paper, we provide results in the \textit{Frames and TV} 
section below,
which 
are an experimental verification of the advantages offered by various frames 
in the CS context. 
More importantly, they show the added benefit of incorporating signal 
structure in 
the sampling procedure, which provides ample freedom to choose the sparsifying 
system for a CS recovery.

\section{A new CS approach in Fluorescence Microscopy}
\label{s:fm}

\def\ppsf{\ensuremath{P_\mathrm{psf}}}
\def\Psim{\ensuremath{\Psi_\mathrm{m}}}
\def\nsq#1{#1^{\raisebox{-1pt}{$\scriptstyle{N}^{\raisebox{-1pt}{$\scriptscriptstyle2$}}$}}}
\def\nnsq#1{#1^{\raisebox{-1pt}{$\scriptstyle\res{N^{\raisebox{-1pt}{$\scriptscriptstyle2$}}}$}}}
\def\allones{\ensuremath{\mathbf{\underline{1}}}}

We start with a FM example which encompasses many of the previously 
enumerated effects and benefits and shows how they allow to improve 
performance in a CS setting. The subsequent sections provide focused 
detailed discussions.

Compressive Fluorescence Microscopy (CFM), a Type~II 
problem, where the sampling operator can be chosen, was first introduced {and 
implemented practically} by 
Studer et al.~\cite{Candes_PNAS} and we refer the reader to their work for 
details. In short, a digital micromirror device of 
$\res{N}$ mirrors can form any $\res{N}$ pattern of 0's and 1's to 
project multiple parallel laser beams (the 1's) through a lens onto a 
specimen, which is excited and emits light (the fluorescence) of a different 
wavelength, collected and summed by a photodetector, thus taking one CS 
measurement. 
Successive measurements are taken, changing the $\res{N}$ pattern on the 
mirrors each time.

\paragraph{Initial approach.}
Studer et al.\ used a 2D Hadamard matrix, i.e.\ 
each pattern on the mirrors corresponds to one row in the 2D Hadamard matrix, 
reordered into $\res{N}$.

\textit{Subsampling pattern.} The pattern $\Omega$ used in~\cite{Candes_PNAS} 
was the ``half-half'' scheme, i.e.\ a two-level sampling scheme where the 
first 
level samples fully and the second level samples uniformly at random.
Hadamard matrices contain $1$'s and $-1$'s but digital mirrors can only 
represent 1's and 0's so~\cite{Candes_PNAS} used the modified sampling 
operator $\Psim=(\Psi+\allones)/2$ where \allones\ is the all-1's matrix and 
$\Psi$ is the Hadamard operator, and solved $\min_z\|\Phi z\|_1$ s.t. 
$\|y-\Psim x\|<\eta$. This, however, is suboptimal since, unlike
$\Psi$, $\Psim$ is non-orthogonal and far from an isometry.

\textit{Point spread effect.} A point emission of light is 
spread by a lens into an airy disc, a blurring effect. The 
lens acts as a circular low-pass filter, i.e.\ its 2D Fourier spectrum is a 
disc.
This is important in the CFM setup where the patterns of 0's and 1's contain 
dicontinuities. To mitigate the lens point-spread effect the above authors 
\textit{binned} mirrors together into groups of \res{2} or \res{4} to 
represent a single 0 or 1 value, which narrows the Fourier response of the 
combined 
light beam coming from such a group, and is less affected by the lens. The 
major drawback is that this reduces the resolution of the recovered image by 2 
or 4 times. This is a serious 
limitation: as we shall see, CS recovery improves with resolution so limiting 
to low resolutions causes a cap in performance. Also, the point-spread 
function (PSF) of the lens was ignored during the CS minimization recovery.

\textit{Photonic noise.} A further challenge is the photonic (Poisson) 
noise 
at the receptor, which essentially counts the number of photons in a preset 
time interval. Unlike white Gaussian noise in other systems, the photonic 
noise mean and variance are signal dependent: the noise power grows with the 
signal. This is accentuated by the CFM setup since there are always $N^2/2$ 
light beams at a time, generating a large background luminance (DC offset) and 
thus impacting the signal-to-noise ratio when measuring higher frequency 
Hadamard components, e.g.\ patterns where 0's and 1's are alternating.

\paragraph{New approach.} In what follows we present an approach which 
employs a multilevel subsampling pattern, explaining why it is beneficial, 
takes into account the lens PSF as well as the photonic noise, 
and also avoids mirror binning, thus reaching much higher resolutions, of the 
order of \res{2048}.
The new approach employs a few techniques to improve performance and 
strives to stay loyal to the planned practical setup in collaboration with the 
Cambridge Advanced Imaging Centre (CAIC), now in the process of building a 
fluorescence microscope of the scale $N=2048$.

Given a subsample pattern $\Omega\subseteq\{1,\ldots,N^2\}$ with $|\Omega|=m$, 
denote with $\ppsf\in\nnsq{\{0,1\}}$ the projection matrix corresponding to 
the Fourier response of the PSF (the disc low-pass filter), with 
$F\in\nnsq{\bbC}$ the 2D Fourier transform and with $x\in\nsq{\bbR}$ the 
original specimen image ordered as a vector. In a noiseless scenario, 
the measurements $\gamma=\{\gamma_i\}\in{\bbN^m}$ would be
\be{\label{e:fm-meas-ideal}
{\gamma=\big\lfloor\big| P_\Omega F^* \ppsf F \Psim x \big|\big\rfloor
	= \left\lfloor\left| \frac{1}{2} P_\Omega F^* \ppsf F 
		\left(\Psi+\allones\right) x \right|\right\rfloor\!\!.}
}
Since the photonic noise is dominant in the CFM setup, the actual 
measurements $y$ can be modelled as values drawn from a Poisson distribution 
with mean and variance equal to $\gamma$, \ie
\be{\label{e:fm-meas-real}
	y=\{y_i\}\sim\text{Poisson}(\{\gamma_i\})\in{\bbN^m}
}
Knowing that a spatially 
wide light beam is very little affected by the PSF, the following 
approximation holds to a high accuracy:
\be{
F^* \ppsf F \,\allones x \;\simeq\; \allones x,
}
since \allones\ gives the widest combined beam, so we can transform the 
measurements $y$ into
\be{
y' = 2y- \underline{\bf y_1},
}
where $\underline{\bf y_1}$ is the vector with all entries equal to $y_1$ 
which 
represents the measurement taken with the all-1's pattern on the mirrors
(the first Hadamard matrix row). This allows us to solve 
\eref{e:l1} using the sampling operator $F^*\ppsf F \Psi$ (or even just 
$\Psi$ as we shall explain), which allows fast transforms and is much 
closer to an isometry compared to 
$F^*\ppsf F\Psim$ (or $\Psim$). Thus we 
solve:
\be{\label{e:l1_cfm}
	\min_{z\in\bbC^{N^2}}\|z\|_1 \quad\text{s.t.}\quad
	\|y'-P_\Omega U z\|<\eta,
}
where $U=F^*\ppsf F \Psi \Phi^*$ and $\Phi$ is an orthobasis like wavelets in 
which the image is expected to be sparse.
This allows us to use fast transforms exclusively for $U$ and its 
adjoint $U^*$, needed during the above minimization, since both 
\ppsf\ and $\Psi$ are real and symmetric, hence self-adjoint, so 
$U^*=\Phi\Psi^* F^*\ppsf F$.

\begin{figure}
\centering
\colwidth{0.315\linewidth}
\def\arraystretch{0}
\newcommand{\abrside}[1]{\vertsidetext{0.15\linewidth}{#1$\times$#1}}
\begin{tabular}{@{}\colpt{4}\colpt{2}\colpt{2}\colpt{0}}
  & \sidetext{Raster scan} & \sidetext{New CS approach} & 
  \sidetext{Initial CS approach}\\[4pt]
  \abrside{256}
  &\overerr{0256_6p3_had_rwtdwtDB4_5_4KDU_ref}
  &\overerr[]{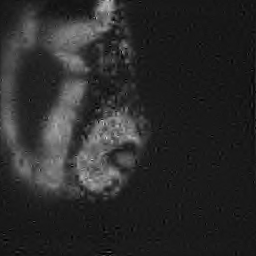}
  &\overerr[]{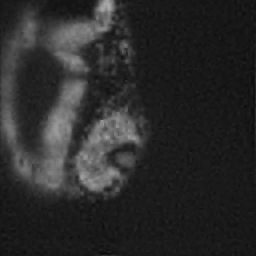}
 \\[3pt]
  \abrside{512}
  &\overerr{0512_6p3_had_rwtdwtDB4_6_A4M9_ref-crop}
  &\overerr[]{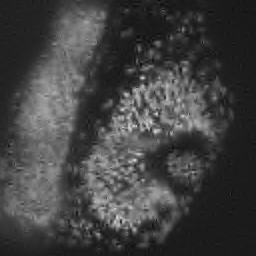}
  &\overerr[]{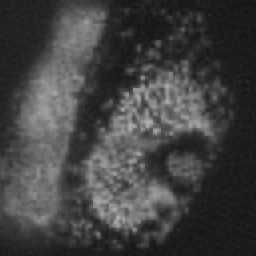}
 \\[3pt]
  \abrside{1024}
  &\overerr{1024_6p3_handle_rwtdwtDB4_7_P25O_ref-crop}
  &\overerr[]{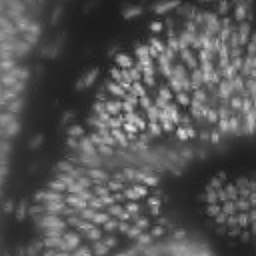}
  &\overerr[]{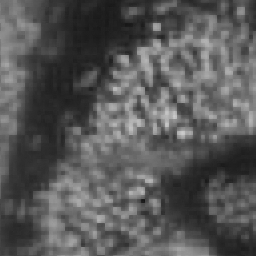}
 \\[3pt]
  \abrside{2048}
  &\overerr{2048_6p3_handle_rwtdwtDB4_8_DKY5_ref-crop}
  &\overerr[]{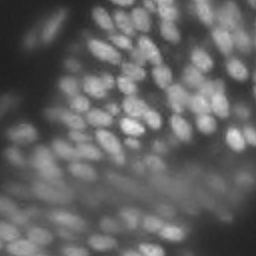}
  &\overerr[]{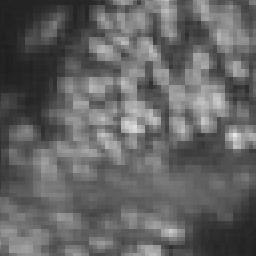}\\
\end{tabular}
\caption{{\capfb Fluorescence Microscopy (FM) example.} Recovering from 6.25\% Hadamard
coefficients \eref{e:fm-meas-real} into Daubechies-4.
The subsampling pattern for the new CS approach is the the one used in 
\fref{f:flip-tests} for 
the FM case.\label{f:fm-full}}
\end{figure}

\fref{f:fm-full} shows zebra fish cells captured by a \res{2048} raster 
scan florescence microscope from CAIC. Measurements $y$ were formed using 
\eref{e:fm-meas-real} (adding large Poisson noise with mean $\gamma$ 
\eref{e:fm-meas-ideal}), where the cutoff 
value of \ppsf\ corresponds to the PSF of the lenses that will be used in 
the microscope which is now being built, so the real-world measurements are 
expected to be close to $y$. As can be seen the 
new approach CS recovery is much improved compared to the initial 
approach that uses $\Psim$ and \textit{half-half} sampling.

\paragraph{Why multilevel subsampling works better.} While in 
\fref{f:fm-full} the use of $F\ppsf F^*\Psi$ and $y'$ (instead of $\Psi_m$ and 
$y$) provide a substantial advantage, a 
multilevel subsampling pattern $\Omega$ intrinsically mitigates the effects of 
the PSF, thus avoiding any mirror binning. This in turn achieves higher 
resolutions and also higher quality recovery. The key, as 
discussed in 
\AAA{the section \textit{New CS principles}}, 
is that $\Omega$ should 
follow the asymptotic decrease of the coherence of Hadamard with wavelets. 
A Hadamard matrix, much like Fourier, captures spatial frequency of increasing 
orders. 
The Hadamard rows that give higher coherence with wavelets correspond to 
lower spatial frequencies, as seen in \fref{f:incoh}, i.e.\ the rows 
with more adjacent 0's or 1's. The new theory states that favoring those rows 
when sampling will provide better CS recovery. Importantly, in the CFM setup 
those rows also inherently emulate mirror binning, owing to the adjaceny of 
1's and 0's, so they are bound to be less affected by the lens PSF. For this 
reason, for an appropriate choice of the multilevel pattern, one could even 
simply use the faster $\Psi$ (instead of the full $F^*\ppsf F\Psi$) as the 
sampling operator. In 
contrast, a ``half-half'' subsampling pattern, besides bound to perform more 
poorly as it does not closely follow the asymptotic incoherence, also 
subsamples heavily from the high spatial frequency rows which are more 
severely affected by the PSF, further decreasing CS recovery quality.

\section{Resolution dependency}
\label{s:mri}

\begin{figure}
\centering
\def\arraystretch{0}
\colwidth{0.31\linewidth}
\newcommand{\abrside}[1]{\vertsidetext{34pt}{#1$\times$#1}}
\begin{tabular}{@{}\colpt{4}\colpt{2}\colpt{2}\colpt{0}}
  & \sidetext{Full sampled} & \sidetext{5\% subsampled} & \sidetext{Subsample 
  map}
  \\[4pt]\abrside{256}
  &\overerr{0256_6p3_dft_rwtdwtDB4_5_peek100_Y1B7_ref-crop}
  &\overerr[10.8]{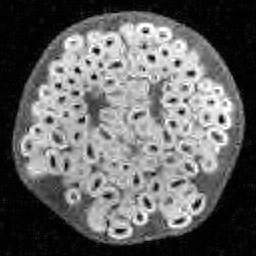}
  &\overerr{0256_6p3_dft_rwtdwtDB4_5_peek100_Y1B7_map}
  \\[3pt]\abrside{512}
  &\overerr{0512_6p3_dft_rwtdwtDB4_6_peek100_BHI2_ref-crop}
  &\overerr[6.01]{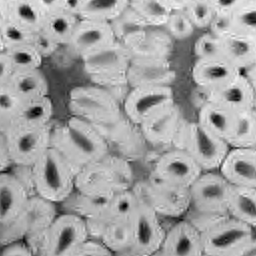}
  &\overerr{0512_6p3_dft_rwtdwtDB4_6_peek100_BHI2_map}
  \\[3pt]\abrside{1024}
  &\overerr{1024_6p3_dft_rwtdwtDB4_7_peek100_KOCM_ref-crop}
  &\overerr[3.60]{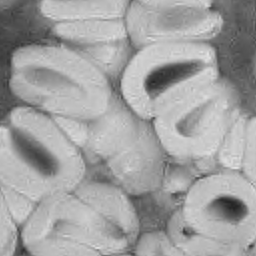}
  &\overerr{1024_6p3_dft_rwtdwtDB4_7_peek100_KOCM_map}
  \\[3pt]\abrside{2048}
  &\overerr{2048_6p3_dft_rwtdwtDB4_8_peek70_SXEG_ref-crop}
  &\overerr[1.87]{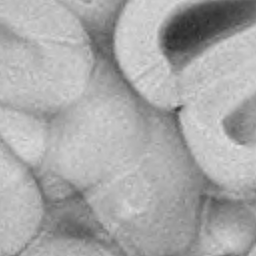}
  &\overerr{2048_6p3_dft_rwtdwtDB4_8_peek70_SXEG_map}
\end{tabular}
\caption{{\capfb MRI example.} Multilevel subsampling of 5\% Fourier coefficients 
recovered into Daubechies-4.\label{f:mri-fraction}}
\end{figure}

\AAA{A resolution dependency effect} could first be 
noticed in \fref{f:fm-full} where the CS recovery gets better as the 
resolution increases since the image is increasingly (asymptotically) sparser 
in wavelets, and the coherence between Hadamard and wavelets decreases 
asymptotically (see \fref{f:incoh}).
\fref{f:mri-fraction} shows a \res{2048} MRI image of a 
pomegranate fruit obtained using a 3T Philips MRI machine, which also contains 
noise specific to MRI. {MRI is an example of Type~I problem, where the 
sampling operator is imposed, but the same new CS principles apply.} We 
subsampled a fixed fraction of 6.25\% Fourier samples and solved 
\eref{e:l1} with $U=\Phi_{\mathrm{dft}}\Psi^*_{\mathrm{dwt}}$. The 
asymptotic sparsity of the wavelet coefficients and the asymptotic incoherence 
of Fourier and wavelets (see \fref{f:incoh}) exploited via multilevel sampling 
again yield increasingly better reconstruction quality as the resolution 
increases, in this case to the point where differences are hardly noticeable. 

\begin{figure}
{\color{black}
\centering
\def\arraystretch{0}
\def\canureadmeh{38pt}
\setlength{\fboxsep}{0pt}
\setlength{\fboxrule}{0pt}
\newcommand{\abrside}[2][]{\raisebox{21pt}{\begin{minipage}[t]{16pt}
	\sidetext{\flushright(#2)\\[2pt]\captext#1}\end{minipage}}}
\newcommand{\abrimg}[1]{\fbox{\includegraphics[height=\canureadmeh,trim={0pt
 	45pt 90pt 15pt},clip]{#1}}}
\newcommand{\abrmap}[1]{\fbox{\includegraphics[height=\canureadmeh]{#1}}}
\begin{tabular}{@{}\colpt{3}\colpt{1.5}\colpt{1.5}\colper{0}}
  & \sidetext{No noise} & \sidetext{White noise (16 dB SNR)} & \sidetext{Map}
  \\[4pt]\abrside[512]{a}
  &\abrimg{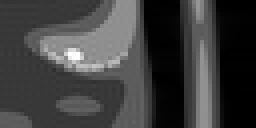}
  &\abrimg{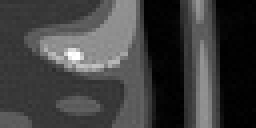}
  &\abrmap{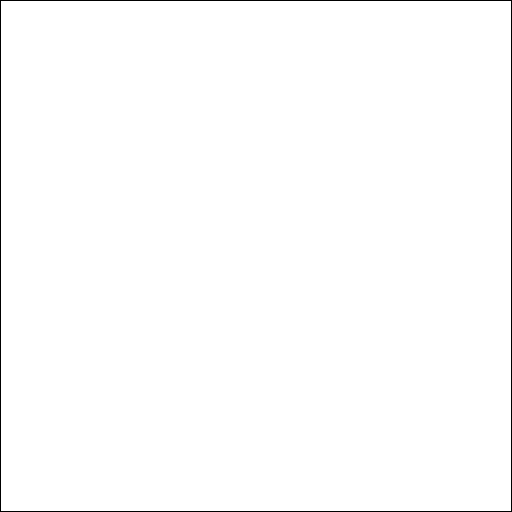}
  \\[3pt]\abrside[2048]{b}
  &\abrimg{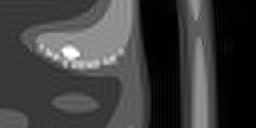}
  &\abrimg{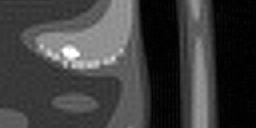}
  &\abrmap{{2048_512_zeropad_map}.png}
  \\[3pt]\abrside[2048]{c}
  &\abrimg{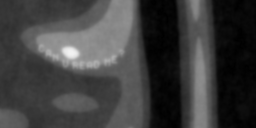}
  &\abrimg{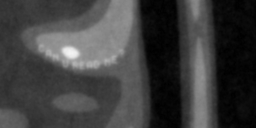}
  &\abrmap{{2048_6.25_db4_61D4_map}.png}
  \\[3pt]\abrside[2048]{d}
  &\abrimg{{2048_6.25_db4_61D4_rec-crop}.png}
  &\abrimg{{2048_6.25_db4_2A32_rec-crop-noise}.png}
  &\abrmap{{2048_6.25_db4_61D4_map}.png}
\end{tabular}
\caption{{\capfb MRI example.} Recovery from a fixed number of 
512\textsuperscript{2} = 262144 Fourier coefficients of the phantom from 
\fref{f:CS_Evol}.
{\bf (a)} \res{512} linear, from the first \res{512} 
coefficients.
{\bf (b)} \res{2048} CS into Daubechies-4,
from the first \res{512} (6.25\%) coefficients.
{\bf (c)} \res{2048} linear, from 512\textsuperscript{2} (6.25\%) 
coefficients sampled with a multilevel sampling map.
{\bf (d)} \res{2048} CS into Daubechies-4, from the same 
512\textsuperscript{2} (6.25\%) coefficients from 
(c).\label{f:glpu-canureadme}}}
\end{figure}

\paragraph{Fixed number of samples.} A more {striking} result of 
asymptotic 
sparsity and asymptotic incoherence is obtained by fixing the 
number of samples taken, instead of the fraction. This was done in 
\fref{f:glpu-canureadme}, sampling the same number of 512\textsuperscript{2} = 
262144 Fourier coefficients in four scenarios, the latter revealing previously 
hidden details when using multilevel sampling from a broader spectrum.

The explanation? The higher resolution opens up higher-order regions of 
wavelet coefficients which are mostly sparse, and higher-order regions of 
coherence between sinusoids and wavelets (see \fref{f:incoh}) which is low. 
As discussed in 
\AAA{the section \textit{New CS Principles},}
when using a 
nonlinear recovery, these two asymptotic effects can be fruitfully 
exploited with a multilevel sampling scheme that spreads the same number of 
samples across a wider range, aiming for the more coherent regions and 
reconstructing finer details to a much clearer extent, even in the presence of 
noise in this example. It is worth noting that other sampling 
strategies (e.g.\ half-half) will also benefit from sampling at higher 
resolutions, provided samples are sufficiently spread out, but a 
multilevel sampling strategy will provide near optimal guarantees 
\cite{AHPRBreaking}.

\begin{figure}[t]
{\color{black}
\centering
\frutiger at 6pt
\colwidth{0.31\linewidth}
\imgstroke{0.25pt}
\begin{tabular}{@{}\colper{0.015}\colper{0.015}\colper{0}}
\overerr{0256_12p6_had_rwtdwtDB4_5_P9LY_ref}&
\overerr{0256_12p6_had_rwtdwtDB4_5_P9LY_map}&
\overerr{2048_12p5_had_rwtdwtDB4_8_ZMN5_map}\\[-2pt]
\captext{Original} & \captext{Best 256$\times$256 map} & \captext{Best 
2048$\times$2048 map}
\end{tabular}
{\includegraphics[height=95pt]{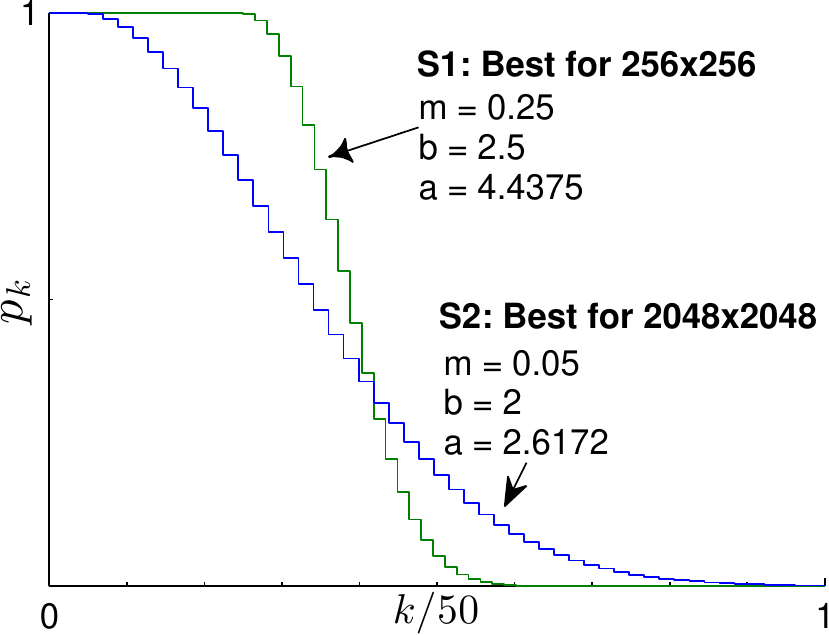}}\quad
\raisebox{7pt}{{\includegraphics[height=88pt]{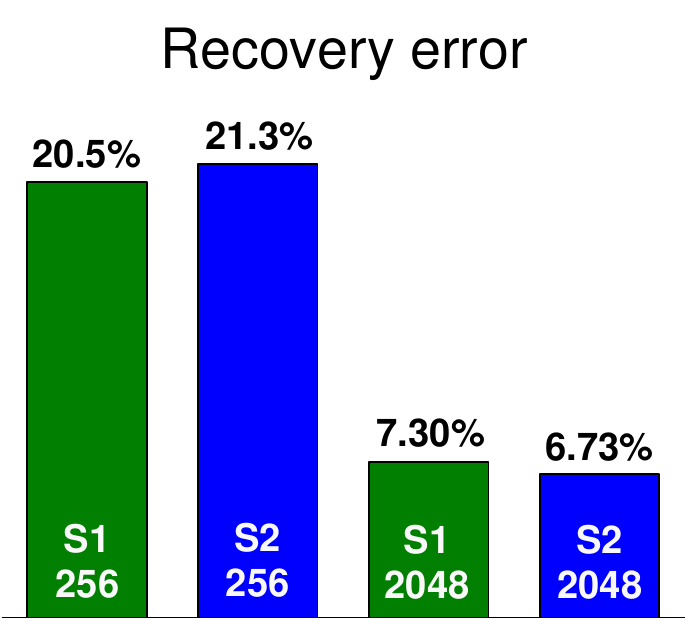}}}\\[5pt]
\caption{{\capfb Compressive Imaging (CI) example.} Recovering from 12.5\% 
Hadamard coefficients into Daubechies-4. The best subsampling maps 
\eref{e:sublaw-gengauss} were 
computed heuristically for each resolution of this 
image.\label{f:res-dep-submap}}}
\end{figure}

\paragraph{{Sampling strategy is also resolution dependent.}} {The 
resolution dependency effect also influences the optimum subsampling 
strategy, in that the latter will depend on the resolution in addition 
to signal sparsity and structure. \fref{f:res-dep-submap} shows an experiment 
in which we computed the best subsampling patterns \eref{e:sublaw-gengauss} 
for two resolutions of the same image. As is evident, the resulting patterns 
are different and also yield different results for the same resolution.}

\paragraph{{Remarks.}} {There are a few remarks worth making 
at this point.\\
~~~-- By simply going higher in resolution (e.g.\ further in the 
Fourier spectrum), one can recover a signal much closer to the exact one, yet 
taking the same number of measurements; or\\
~~~-- By simply going higher in resolution one can obtain 
the same reconstruction quality, yet taking {fewer} measurements.\\
~~~-- These experiments showed that it is important in practice to be 
able to access 
high resolutions (higher frequencies in the MRI case) in order for CS to 
provide higher gains. Thus, for MRI, the benefits will be even more visible on 
future generations of MRI machines.\\
~~~-- Multilevel sampling can better exploit the resolution dependency effect 
and allows for better tailoring according to sparsity 
structure, resolution or application specific requirements (e.g.\ different 
patterns for different body parts, allow lower 
overall quality but recover contours better etc) as opposed to uniform random 
sampling or sampling schemes such as half-half~\cite{Candes_PNAS} or 
continuous power laws~\cite{KrahmerWardCSImaging}.\\
~~~-- Thirdly, practical CS in MRI has several limitation regarding
point-wise sampling. The multilevel patterns used here are the result of our 
quest for a theoretically optimal sampling pattern, which could then be 
approximated by realistic MRI patterns or contours, e.g.\ parametric spirals. 
The latter is work in progress in collaboration with the Wolfson Brain Imaging 
Centre.}

\section{\AAA{Infinite dimensional problems}}
\label{s:tem}

{The underlying model in some applications is continuous, such as in MRI, 
EM, X-ray tomography and its variants. These are Type~I problems, where 
the sensing operator is imposed. In MRI,} the measurements $y$ are samples of 
the continuous (integral) Fourier 
transform $\cF$. The same applies for EM and X-ray and its variants, where 
the Radon transform is sampled one angle at the time. Via the Fourier slice 
theorem, the procedure is equivalent to sampling the Fourier transform at 
radial lines and so the Fourier and Radon transform recovery problems are 
equivalent to recovering the continuous $f$ from pointwise samples {$\hat g$, 
which are evaluations of}
\be{
g = \mathcal{F}f, \quad f \in \rL^2(\mathbb{R}^d), \quad \mathrm{supp}(f) 
\subseteq [0,1]^d.
\label{e:invfourier}
}
\paragraph{The issue.} 
Consider the CS recovery into wavelets (see 
\cite{Adcock_Hansen_Book} for the general case). Using discrete tools, 
in this case $U=\Psi_\mathrm{dft}\Phi^*_\mathrm{dwt}$ in \eref{e:l1}, to 
solve a continuous inverse problem can introduce notable errors, due to
\textit{measurement mismatch} and the \textit{wavelet 
crime}~\cite{StrangNguyen}.
The former assumes a discrete model for $f$, $\tilde f = \sum_{j=1}^N
\tilde{\beta}_j\psi_j$, where $\psi_j$ are step functions, and 
then (more seriously) replaces the continuous $f$ and $\cF$ in 
\eref{e:invfourier} with {$\tilde f$} and DFT
respectively, leading to the discretization $\tilde g=\Psi_\mathrm{dft}
\tilde{\beta}$, which is a poor approximation {of the samples of $g$}. The 
wavelet crime is as follows. Given scaling and mother wavelet functions 
$\varphi$ and $\psi$, obtaining the wavelet coefficients of $f$ via the DWT 
should assume $f=\sum_{j=-\infty}^\infty \beta_j\varphi(\cdot-j)$ and then 
compute them from $\{\beta_j\}$ via the DWT. The crime is when one simply 
replaces $\beta_j$ with pointwise samples of $f$.

{To illustrate these issues, let $y\!=\!P_\Omega\hat{g}$} be the measurements 
in 
\eref{e:l1}, where {$\hat{g}$} are the first $2N$ 
continuous Fourier samples of 
$f$. The matrix $\Psi^*_\mathrm{dft}$ maps {$\hat{g}$} to a vector 
$x\in\bbC^{2N}$ on an equispaced $2N$ grid of points in $[0,1]$. Specifically, 
$\Psi_\mathrm{dft}x={\hat{g}}$ 
where $x$ is given by the values $f_N(t)\!=\!1/2 \sum^{N}_{j=1-N} \cF f(j/2) 
\E^{2 \pi \I \epsilon j} $ on the $2N$ grid. Taking $x_0=\Phi_\mathrm{dwt}x$, 
the right-hand side of \eref{e:l1} 
becomes $P_{\Omega}\Psi_\mathrm{dft}\Phi^*_\mathrm{dwt}x_0$. But \eref{e:l1} 
now 
requires $x_0$ to be sparse, which means the truncated 
Fourier series $f_N$ must be sparse in wavelets, which cannot happen. While 
$f$ is assumed sparse in wavelets, $f_N$ is not, since it 
consists of smooth complex exponentials. {Large errors thus occur in the 
recovery, as there is no sparse solution due to the poor approximation of $f$ 
by $f_N$. This could be avoided if the measurements arose from the DFT, but 
that would be the well-known and pervasive \textit{inverse crime}~\cite{GLPU}: 
artificially superior performance when data is simulated incorrectly using the 
DFT, as opposed to the continuous FT, which is the true underlying model.}

\paragraph{How to solve.} 

\begin{figure}
\centering
\capfontsize{7pt}
\newcommand{\abrdesc}[3]{\captext{#1}&\captext{#2}&\captext{#3}}
\colwidth{0.31\linewidth}
\imgstroke{0pt}
\begin{tabular}{@{}\colper{0.01}\colper{0.01}\colper{0}}
  \overerr{006p149_map}
  &\overerr{006p149_orig}
  &\overerr{006p149_orig_zoom}\\[-1pt]
  \abrdesc{Sampling map}{Original}{Enlarged}\\[5pt]%
   \overerr[17.8]{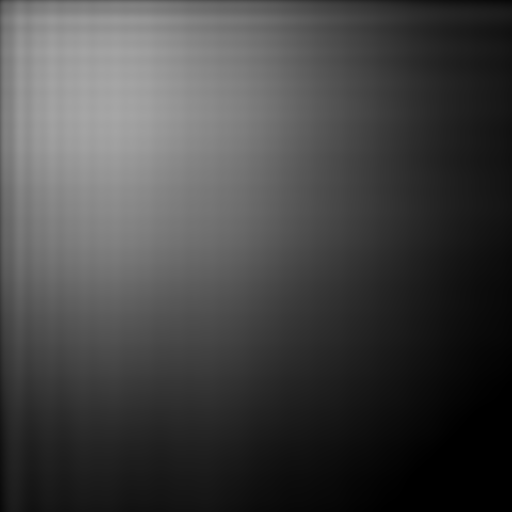}
  &\overerr[9.8]{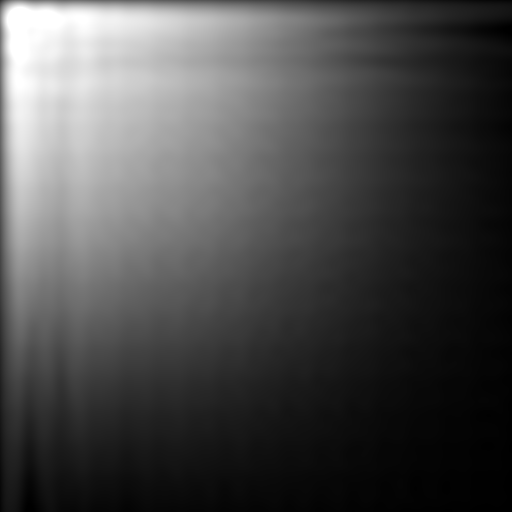}
  &\overerr[0.1]{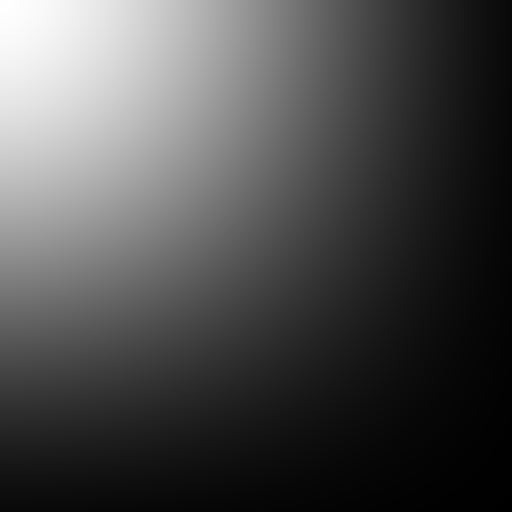}\\[-1pt]
  \abrdesc{Linear, Inverse DFT}{CS $\ell^1$, DFT to}{CS $\ell^1$, Inf-dim 
  CS,}\\[-2pt]
  \abrdesc{}{Periodic DB6}{Boundary DB6}
\end{tabular}
\caption{{\capfb Electron Microscopy (EM) example.} Recovery from 16120 (6.15\%) 
continuous FT samples onto a \res{512} grid.\label{f:infdim}}
\end{figure}

The above crimes stem from discretizing first and 
then applying finite-dimensional tools. Instead, we shall use the techniques 
of \textit{infinite-dimensional} CS~\cite{BAACHGSCS}, i.e.\ first formulate 
the 
problem in infinite dimensions and \textit{then} discretize.
Let $\{ \psi_j \}_{j \in \bbN}$ and $\{ \varphi_j \}_{j \in \bbN}$ be the 
sampling and sparsity bases (Fourier and wavelets). If $f = \sum_{j \in 
\bbN} \beta_j \varphi_j$, 
then the unknown vector of 
coefficients $\beta = \{ \beta_j \}_{j \in \bbN}$ is the solution of $U \beta 
= \hat{f}$, where
\bes{
U =
  \left(\begin{array}{ccc} \left < \varphi_1 , \psi_1 \right >   & \left < 
  \varphi_2 , \psi_1 \right > & \cdots \\
\left < \varphi_1 , \psi_2 \right >   & \left < \varphi_2 , \psi_2\right > & 
\cdots \\[-4pt]
\vdots  & \vdots  & \ddots   \end{array}\right),
}
and $\hat{f} = \{ \hat{f}_j \}_{j \in \bbN}$ is the infinite vector of 
samples of $f$. Suppose we have access to a finite number of samples
$\{ \hat{f}_j : j \in \Omega \}$, where $\Omega$ is the sampling map.
To recover $\beta$ from these samples, we first formulate the 
infinite-dimensional optimization problem
\be{
\label{e:l1-inf}
\inf_{z \in \ell^1(\bbN)} \| z \|_1 \text{ subject to }
P_{\Omega} Uz = P_{\Omega} \hat{f},
}
which commits no crimes. However, \eref{e:l1-inf} cannot typically be solved 
numerically, so we now discretize. We
truncate $U$ to $K\in\bbN$ columns, solving the now finite-dimensional 
problem
\be{
\label{findimopt}
\min_{z \in P_K(\ell^2(\bbN))} \| z \|_1 \text{ subject to }
P_{\Omega} U P_K z = P_{\Omega} \hat{f}.
}
We refer to this as infinite-dimensional CS, where $K$ should be sufficiently 
large for good recovery. \fref{f:infdim} shows such 
a problem where the continuous function 
$f(x,y)=\exp(-x-y)\cos^2\!\!\;(x)\cos^2\!\!\;(y)$ is recovered on a \res{512} 
grid 
from 
16120 (6.15\%) continuous Fourier samples taken radially, as in an EM setup. 
The infinite-dimensional CS recovery is done in boundary DB6 wavelets instead 
of the {periodic} DB6 (to preserve the vanishing moments at the boundaries). 
It 
is evident that its 
reconstruction is far superior to both the discrete linear reconstruction (via 
$\Psi^*_\mathrm{dft}$) and discrete CS reconstruction (via 
$U=\Psi_\mathrm{dft}\Phi^*_\mathrm{dwt}$) which are affected by the crimes.

In conclusion, given sufficiently many vanishing moments, the infinite 
dimensional CS with boundary wavelets will give substantially better 
convergence than the slow truncated Fourier convergence $f_N$ whenever $f$ is 
non-periodic. Knowing that the computational complexity in this case is the 
same as with FFT, this means that infinite-dimensional CS yields a markedly 
better approximation of $f$ at little additional cost. The 
infinite-dimensional CS approach described here is of particular benefit to 
applications like EM and spectroscopy where smooth functions are encountered. 
{This is work in progress in collaboration with the Department of Materials 
Science \& Metallurgy at University of Cambridge.}

\section{Structured sampling vs Structured recovery}
\label{s:algs}

The new CS principles in this paper take into account sparsity 
structure in the sampling procedure via multilevel sampling of non-universal 
sensing matrices.
Sparsity structure can also be taken into account in the recovery algorithm. 
An example of such an approach is model-based CS
\cite{BaraniukModelCS}, which assumes the signal is piecewise smooth and 
exploits the connected tree structure of its wavelet coefficients to reduce 
the search space of the matching pursuit algorithm~\cite{cosamp}. Another 
approach is the class of message passing and approximate message passing 
algorithms~\cite{BaraniukMessagePassing,DonohoAMP}, which exploit the 
persistence across scales structure~\cite{mallat09wavelet} of wavelet 
coefficients by a modification to iterative thresholding algorithms inspired 
by belief propagation from graph models. This can be coupled with hidden 
Markov trees to model the wavelet structure, {such as in the Bayesian 
CS~\cite{HeCarinStructCS,HeCarinTreeCS} and TurboAMP~\cite{TurboAMP} 
algorithms.} \BBB{Another approach is to assign non-uniform weights to the 
sparsity coefficients \cite{HassibiWeightedL1}, to favor the important 
coefficients during $\ell^1$ recovery by assuming some typical decay rate of 
the coefficients.}
{A recent approach assumes the actual signal (not its representation 
in a sparsity basis) is sparse and random,
and shows promising theoretical results when using spatially coupled matrices 
\cite{WuVerduLimitsAnalogCompression,KrzakalaPhysRev,DonohoSpatialCouplingAMP},
 yet it is unclear how a real-world setup can be implemented 
 where signals are sparse in a transform domain.}

The main difference is that the former approach, i.e.\ multilevel sampling of 
asymptotically incoherent matrices, incorporates sparsity structure in the 
sampling strategy and uses standard $\ell^1$ minimization algorithms, whereas 
the latter approaches exploit structure by modifying the recovery algorithm 
and use universal sampling operators which yield uniform incoherence (see 
\AAA{section \textit{Structure vs Universality}),}
e.g.\ random Gaussian or Bernoulli.

\paragraph{Comparison.} \textit{Structured recovery:}
Due to the usage of universal operators and assumptions on the
sparsity basis, this approach is typically restricted to Type~II 
problems, where the sensing operator can be designed, and is further 
restricted in the choice of the sparsity frame, whose structure is 
exploited by the modified algorithm.

\noindent\textit{Structured sampling:}
In contrast, this approach practically has no limitation regarding the 
sparsity frame, thus allowing for further improvement of CS recovery (see
\AAA{section \textit{Frames and TV}),} 
and it also works for Type~I problems, where the sensing 
operator is imposed.

\begin{figure}
\colwidth{0.45\linewidth}
\capfontsize{7pt}
\newcommand{\abrdesc}[2]{\captext{#1}&\captext{#2}\\[5pt]}
\centering
\begin{tabular}{@{}\colper{0.05}\colper{0}@{}}
\overerr{0256_12p6_had_curve_4_NH5Y_ref}&
\overerr[16.0]{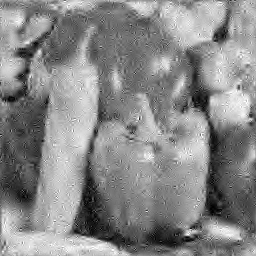}\\[-2pt]
\abrdesc{Original}{Random Bernoulli to db4 --- $\ell^1$}
\overerr[21.2]{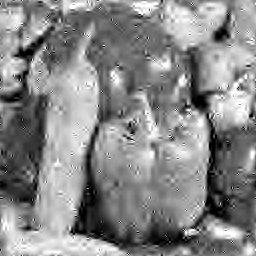}&
\overerr[17.5]{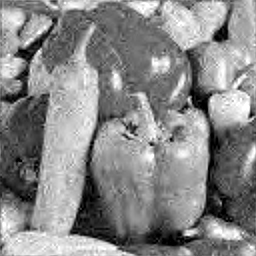}\\[-2pt]
\abrdesc{Rnd. Bernoulli to db4 --- ModelCS}{Rnd. Bernoulli to db4 --- TurboAMP}
\overerr[10.6]{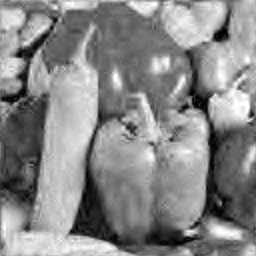}&
\overerr[10.4]{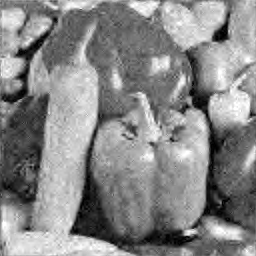}\\[-2pt]
\abrdesc{Random Bernoulli to db4 --- BCS}{\BBB{Rnd. Bernoulli to db4 --- 
Weighted $\ell^1$}}
%
\overerr[7.1]{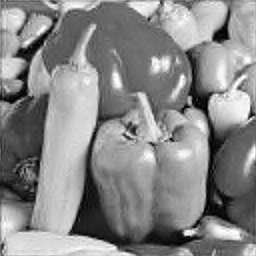}&
\overerr[6.3]{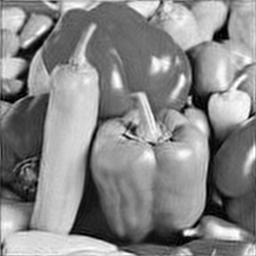}\\[-2pt]
\abrdesc{Multilevel Hadamard to db4 --- $\ell^1$}{Multilevel Had.\ to 
Curvelets --- $\ell^1$}%
\end{tabular}%
\caption{{\capfb Compressive Imaging (CI) example.} 12.5\% 
subsampling at \res{256}. The multilevel subsample map is the one from 
\fref{f:flip-tests}. \CCC{The weighted $\ell^1$ weights were $w_k=a^k$ for all coefficients in the $k$th wavelet scale, where $a=2$ was computed heuristically for this example.}\label{f:algs}}
\end{figure}

To compare performance, we ran a large set of simulations of a CI setup. CI 
\cite{SinglePixel,LenslessBellLabs}, a Type~II problem, is an application 
where universal sensing matrices have been traditionally favored. Here the 
measurements $y$ are typically taken using a sensing matrix with only two 
values (usually $1$ and $-1$). Any matrix with only two values fits this 
setup, e.g.\ Hadamard, random Bernoulli, Sum-To-One (see 
\AAA{section \textit{Storage/speed}}),
hence we can directly compare the two approaches.
\fref{f:algs} shows a representative example from that set, which points 
to the conclusion that asymptotic incoherence combined with multilevel 
sampling of highly non-universal sensing matrices (e.g.\ Hadamard, Fourier) 
allows structured sparsity to be better exploited than universal sensing 
matrices, even when structure is accounted for in the recovery algorithm. 
{The figure also shows the added benefit of being able to use a better 
sparsifying system, in this case curvelets.}

\section{Structure vs Universality: Is universality desirable?}
\label{s:universality}

{Universality is a reason for the popularity in traditional CS of random 
sensing matrices, e.g.\ Gaussian or Bernoulli}. A random matrix $A \in 
\mathbb{C}^{m\times N}$ is universal if for any isometry $\Psi \in 
\mathbb{C}^{N\times N}$, the matrix $A\Psi$ satisfies the RIP with high 
probability. For images, a common choice is $\Psi = \Psi^*_{\mathrm{dwt}}$, 
the inverse wavelet transform. {Universality helps when the signal is 
sparse but possesses no further structure.

But is universality desirable in a sensing matrix when the signal is 
structured? First, random matrices are largely inapplicable in Type~I 
problems, where the sampling operator is imposed. They are applicable mostly 
in Type~II problems, where there is freedom to design the sampling operator. 
Should then one use universal matrices there? We argue that universal matrices 
offer little room to exploit extra structure the signal may have, even in 
Type~II problems.

Practical applications typically entail signals that exhibit far more 
structure than sparsity alone: {in particular,} asymptotic sparsity 
structure in some sparsity basis. Thus,} an alternative is to use a 
non-universal sensing matrix, such as Hadamard, $\Phi_{\mathrm{Had}}$. As 
previously discussed and shown in \fref{f:incoh}, 
$U\!=\!\Phi_{\mathrm{Had}}\Psi^*_{\mathrm{dwt}}$ is completely 
coherent with all wavelets yet asymptotically incoherent, and thus perfectly 
suitable for a multilevel sampling scheme which can exploit the inherent 
asymptotic sparsity. This is precisely what we see in \fref{f:algs}: 
multilevel sampling of a Hadamard matrix can markedly outperform solutions 
employing universal matrices {in Type~II problems. For Type~I problems, 
an important practical aspect is that {many} imposed sensing operators 
{happen to be} highly non-universal and asymptotically incoherent with popular 
sparsity bases, and thus easily exploitable using multilevel sampling, as 
seen in \freftwo{f:CS_Evol}{f:flip-tests}.}

\paragraph{Asymptotic incoherence vs Uniform incoherence.}\label{s:incoh}
The reasons for the superior results are rooted in the 
incoherence structure. Universal and close to universal sensing 
matrices typically provide a relatively low {and flat} coherence pattern.
This allows sparsity to be exploited by sampling 
uniformly at random but, by definition, these matrices cannot 
exploit the distinct asymptotic sparsity structure when 
using a typical ($\ell^1$ minimization) CS reconstruction.

In contrast, when the sensing matrix provides a coherence pattern that aligns 
with the signal sparsity pattern, one can fruitfully exploit such 
structure.  As discussed in 
\AAA{the sections \textit{New CS principles} and \textit{Resolution 
dependency}}
a multilevel sampling scheme is likely to give superior results 
by sampling more in the coherent regions, where the signal is also typically 
less sparse. Even though the optimum sampling strategy is signal dependant 
(see
\AAA{section \textit{Traditional CS}}),
real-world signals, particularly images, 
share a {fairly} common structure and thus good, all-round multilevel 
sampling strategies can be designed. {An added benefit of multilevel 
sampling is that it also allows {for} tailoring {of} the sampling pattern to 
target {application-specific} features rather than an all-round approach, 
e.g.\ allow{ing} 
a slightly lower overall quality but recover{ing} contours better.}


\begin{figure}
\colwidth{0.28\linewidth}
\newcommand{\abrside}[1]{\raisebox{40pt}{\begin{minipage}[t]{0.05\linewidth}%
	\flushright\sidetext{to\\#1}\end{minipage}}}%
\centering
\begin{tabular}{@{}\colper{0.02}\colper{0.01}\colper{0.01}\colper{0}}
& \sidetext{Rand Bernoulli} & \sidetext{Rand Gaussian} & 
\sidetext{STOne}\\[1pt]
\abrside{db4}&
\overerr[15.9]{0256_12p5_bern_dwtdwtDB4_8_A9WU_rec}&
\overerr[15.9]{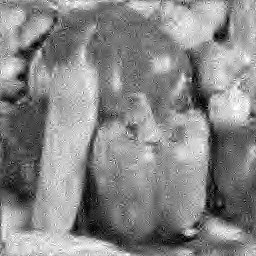}&
\overerr[15.9]{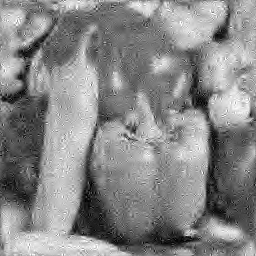}\\
\abrside{db4\\flip}&
\overerr[15.9]{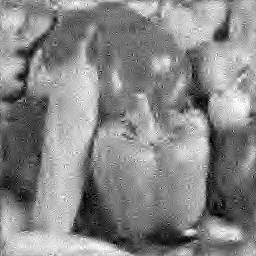}&
\overerr[15.9]{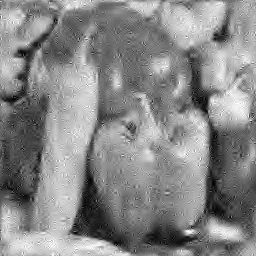}&
\overerr[15.8]{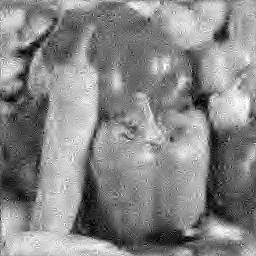}
\end{tabular}%
\caption{Subsampling 12.5\% at \res{256} with 
Bernoulli, Gaussian and STOne; part of a larger set of experiments using 
various images, ratios and sparsity bases. All three matrices yield 
very similar quality, indicating that they behave the same and that 
\textit{universality} and \textit{RIP} hold for all three.\label{f:rand-stone}}
\end{figure}

\section{Storage/speed: Is non-random/orthogonality enough?}
\label{s:storage-speed}

Random matrices have another important practical drawback: they require 
(large) storage and lack fast transforms. This limits the maximum signal 
resolution and yields slow recovery. For example, a \res{1024} experiment with 
25\% subsampling of a random Gaussian matrix would require 2 Terabytes of free 
memory and {$\mathcal{O}(10^{12})$ time complexity,} 
making it impractical at best.


A low maximum resolution is a big limitation not just for
computations. As seen in 
\AAA{the section \textit{Resolution dependency}},
at low resolutions the asymptotic sparsity has not kicked in 
and CS yields only marginal benefits. In order to obtain better 
recovery it is thus of great interest to be able to access high signal 
resolutions, yet random matrices prevent that.

But is the problem of quality CS recovery solved if we {address} the 
storage and speed issues? These were in fact addressed to various 
extents, e.g.\ pseudo-random column permutations of the 
columns of orthogonal matrices such as (block) Hadamard or Fourier
\cite{ScrambledFourier,HadamardBlockScrambled}, Kronecker 
products of random matrix stencils~\cite{BaraniukKronecker}, or even fully 
orthogonal matrices such as the Sum-To-One (STOne) matrix~\cite{STOne} which 
has a fast ${\cal O}(N\log N)$ transform. However, all 
these matrices strive to provide universality, i.e.\ to behave like purely 
random matrices. \fref{f:rand-stone} shows an extract of a large experiment 
\BBB{on various images, resolutions, sparsity bases and subsampling fractions},
which tests for universality and RIP by performing the flip test (see 
\AAA{section \textit{Traditional CS}}).
It is evident that the orthogonal STOne matrix behaves like random matrices in 
the CS context \BBB{(though we note that the STOne matrix was invented to 
serve other purposes as well \cite{STOne}),}
which we probed in \fref{f:algs}. The poor performance is due to their flat 
incoherence with 
sparsity bases, which was discussed in 
\AAA{the section \textit{Structure vs Universality}}.

%


Another solution to the storage and speed problem is to instead use orthogonal 
and structured matrices like Hadamard, DCT or DFT. These have fast 
transforms and do not require storage, but 
also provide asymptotic incoherence with sparsity bases, thus a multilevel 
subsampling scheme can be used. {The important added benefits are that 
this yields significantly better CS recovery in most Type~II problems when 
compared to universal matrices, as discussed in 
\AAA{the section \textit{Structure vs Universality}}
and probed in \fref{f:algs}, and that, unlike universal matrices, it is also 
applicable to Type~I problems, which impose the sensing operator.}

In conclusion, the sensing matrix must contain additional structure 
besides simply being non-random and/or orthogonal%
 in order to provide 
asymptotic incoherence. Typically, sensing and sparsifying 
matrices that are discrete versions of integral transforms, e.g.\ Fourier, 
wavelets \etc will provide asymptotic incoherence, but other orthogonal and 
structured matrices like Hadamard will do so too.

\section{Frames and TV: Freedom of choice}
\label{s:frames}

\begin{figure}
\centering
\colwidth{0.33\linewidth}
\def\arraystretch{0}
\newcommand{\abrover}[2]{\overtext[l]{#1}{#2}}
\begin{tabular}{@{}\colper{0.005}\colper{0.005}\colper{0}}%
\abrover{Subsample map}{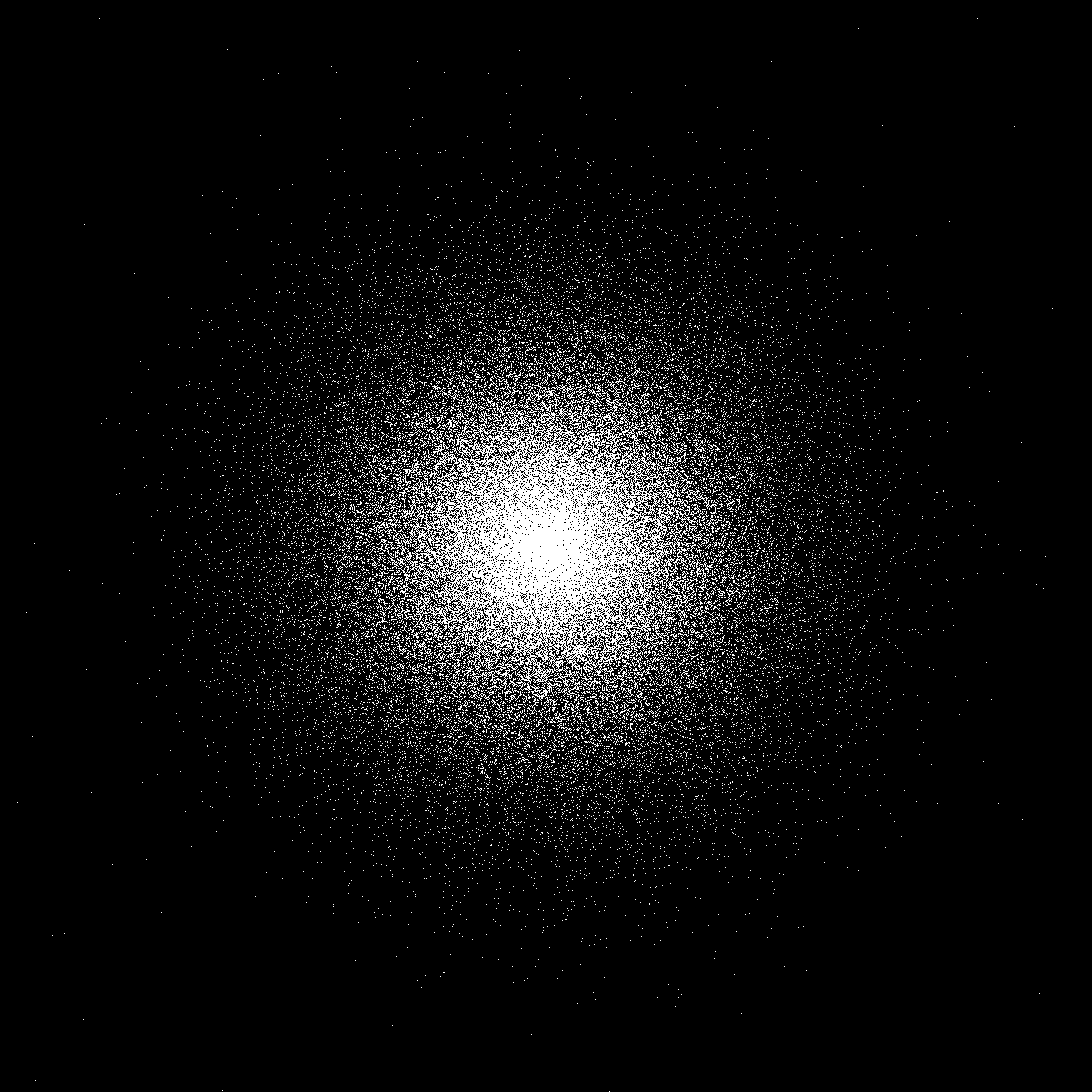}&
\abrover{Original}{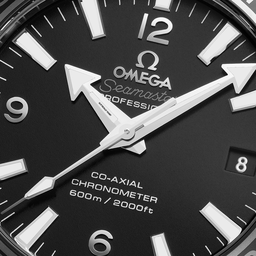}&
\abrover{\vbox{\hbox{\strut%
	Original}\hbox{\strut100\%zoom}}}{2048_watch_full-crop}\\[0.005\linewidth]
\abrover{\vbox{\hbox{\strut Linear}\hbox{\strut inverse 
DFT}}}{2048_6p3_dft_linear-crop}&
\abrover{TV}{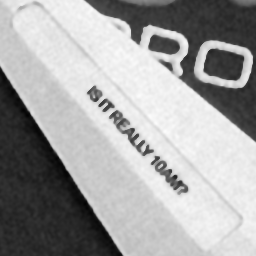}&
\abrover{Daubechies\,4}{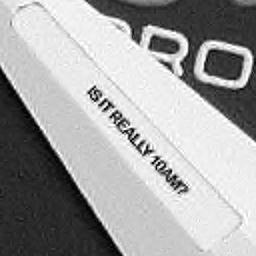}\\[0.005\linewidth]
\abrover{Curvelets}{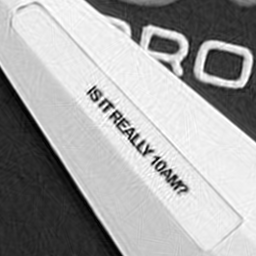}&
\abrover{Contourlets}{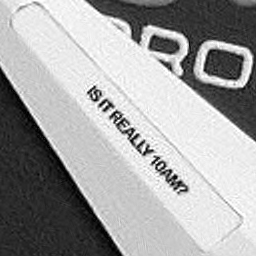}&
\abrover{Shearlets}{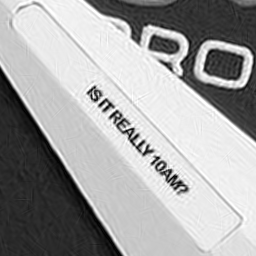}
\end{tabular}
\caption{Recovering from the same 6.25\% DFT coefficients at 
\res{2048}.\label{f:frames}}
\end{figure}

Though it is not the purpose of this paper to investigate frames or TV as 
sparsifying systems in CS, we provide some further 
results. Without going into details, an aspect of interest here
is that many images are known to be relatively sparser in TV and 
frames such as curvelets~\cite{Cand}, shearlets~\cite{Gitta} or contourlets  
\cite{Vetterli}, than in orthobases such as wavelets or DCT.
The results from \fref{f:frames} provide clear experimental 
verifications of the improvements offered by such sparsifying systems 
at practical resolution levels.

More important is that, unlike the class of 
modified recovery algorithms from 
\AAA{the section \textit{Structure sampling vs Structure recovery}},
and in 
addition to the benefits discussed in
\AAA{the sections \textit{Structure vs Universality} and 
\textit{Storage/speed}},
incorporating sparsity structure in the sampling 
procedure also offers complete freedom in the choice of the sparsity 
system. This holds generally, and is of particular interest in applications 
where the sampling operator is imposed.


\section{Concluding remarks}

The traditional CS pillars: sparsity, incoherence and uniform 
random subsampling, are {often} inapplicable in Type~I problems, where the 
sampling operator is imposed (MRI, EM, Tomography, Interferometry etc.), while 
for Type~II problems, where the sampling operator can be designed 
(FM, CI etc.), they provide little room to exploit extra sparsity structure 
that 
real-world {signals} typically possess. This is due to the coherent nature of 
Type~I problems and of the uniform incoherence of universal sampling 
operators with sparsity bases.

The new CS principles: asymptotic 
sparsity, asymptotic incoherence and multilevel subsampling, introduced by the 
authors~\cite{AHPRBreaking} to bridge the gap between theoretical and 
practical CS, provide a better fit for both types of problems. 
This paper shows how the new principles can be used {to} better understand the 
underlying phenomena in practical CS problems, and that an approach based on 
the new CS principles coupled with non-universal sampling operators can 
overcome {many traditional CS limitations} and provide several important 
benefits and improved CS recovery in real-world applications.

%

\begin{acknowledgments}
The authors thank Kevin O'Holleran for providing the 
zebra fish FM image, Andy Ellison for providing the MRI pomegranate image and 
Clarice Poon for the original generalised sampling code. BR and AH acknowledge 
support from the UK Engineering and Physical Sciences Research Council (EPSRC) 
grant EP/L003457/1. AH acknowledges support from a Royal Society University 
Research Fellowship. BA acknowledges support from the NSF DMS grant 1318894.
\end{acknowledgments}

\bibliographystyle{pnas}
\bibliography{MoreForLess}

\end{article}

\end{document}